\documentclass[5p]{elsarticle}

\usepackage{lineno,hyperref}
\usepackage{amsthm,amsmath,amssymb}
\usepackage{mathrsfs}
\usepackage{amsmath}
\usepackage{bm}
\usepackage{color}
\usepackage{subfigure}
\usepackage{multirow}
\usepackage{rotating}
\usepackage{url}
\usepackage{subeqnarray}
\usepackage{cases}
\usepackage{amsmath}
\usepackage{nomencl}
\usepackage{ifthen}
\usepackage[ruled,vlined]{algorithm2e}

\renewcommand{\nomgroup}[1]{%
	\ifthenelse{\equal{#1}{D}}{\item[\textbf{Sets/Indices}]}{%
		\ifthenelse{\equal{#1}{P}}{\item[\textbf{Parameters}]}{%
			\ifthenelse{\equal{#1}{V}}{\item[\textbf{Variables}]}{}}}
}
\newcommand{\figref}[1]{Fig.~\ref{#1}}

\makenomenclature

\journal{Journal of Energy Conversion and Management}

\begin{document}

\begin{frontmatter}

\title{A two-stage robust optimization approach for oxygen flexible distribution under uncertainty in iron and steel plants}

\author[1,2]{Sheng-Long Jiang\corref{cor1}}%
\ead{sh.l.jiang.jx@gmail.com}
\author[3]{Gongzhuang Peng}
\ead{gzpeng@ustb.edu.cn}
\author[2]{I. David L. Bogle}
\ead{d.bogle@ucl.ac.uk}

\cortext[cor1]{Corresponding author: +86 13520412520}
\address[1] {College of Materials Science and Engineering, Chongqing University, Chongqing, 400044, P.R. China}
\address[2] {Center for Process Systems Engineering, Department of Chemical Engineering, University College London, London  WC1E 7JE, U.K.}
\address[3]  {Engineering Research Institute, University of Science and Technology Beijing, Beijing 100083, P.R. China}

\begin{abstract}
Oxygen optimal distribution is one of the most important energy management problems in the modern iron \& steel industry. Normally, the supply of the energy generation system is determined by the energy demand of manufacturing processes. However, the balance between supply and demand fluctuates frequently due to the uncertainty arising in manufacturing processes. In this paper, we developed an oxygen optimal distribution model considering uncertain demands and proposed a two-stage robust optimization (TSRO) with a budget-based uncertainty set that protects the initial distribution decisions with low conservatism. The main goal of the TSRO model is to make “wait-and-see” decisions maximizing production profits and make “here-and-now” decisions minimizing operational stability and surplus/shortage penalty. To represent the uncertainty set of energy demands, we developed a Gaussian process (GP)-based time series model to forecast the energy demands of continuous processes and a capacity-constrained scheduling model to generate multi-scenario energy demands of discrete processes. We carried out extensive computational studies on TSRO and its components using well-synthetic instances from historical data. The results of model validation and analysis are promising and demonstrate our approach is adapted to solve industrial cases under uncertainty.
\end{abstract}


\end{frontmatter}


\section{Introduction}

Energy plays a vital role in the plant-wide optimization of manufacturing systems in the context of achieving a carbon-neutral economy. However, there are two major challenges exposed to most energy-intensive manufacturers: (1) rising energy costs, (2) strict requirement on emission. Especially, these challenges to the iron \& steel industry are even more severe because it is the second-largest energy consumer among the global industrial sectors \cite{johansson2016effects,sun2020material}. Oxygen is one of the most fundamental energy resources for a iron \& steel plant due to some processes need to consume a large amount of oxygen gas, e.g., steelmaking by basic oxygen furnace (BOF), ironmaking by blast furnace (BF). Moreover, oxygen is also a critical health resource, especially in the COVID-19 Pandemic. Therefore, how to make optimal decisions on oxygen production and distribution significantly benefits productivity, energy-efficiency and even indirectly contributes to mitigate public health crisis.

\begin{figure*}[htbp!]
	\centering
	\includegraphics[width=0.85\textwidth]{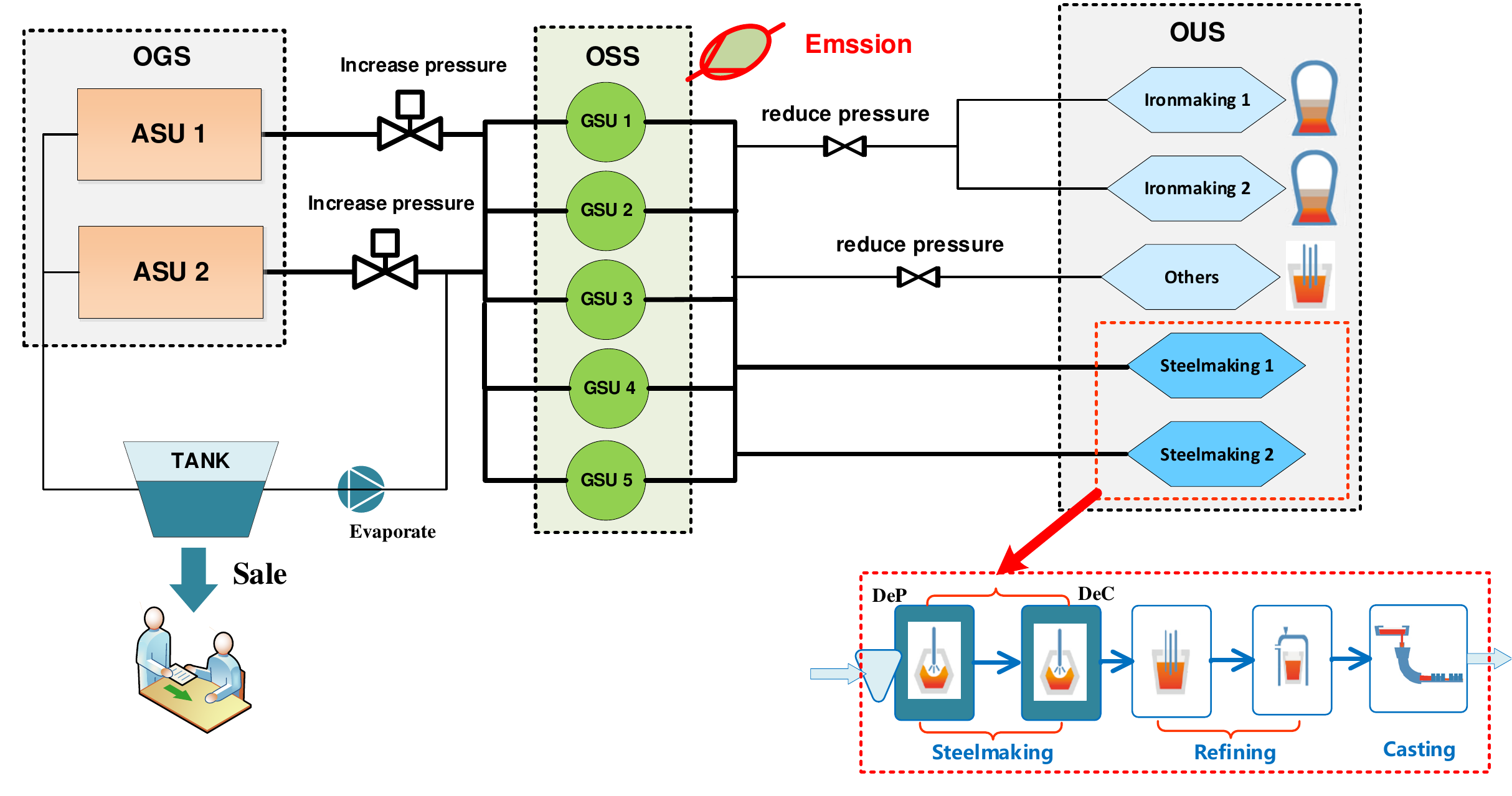}
	\caption{Oxygen supply and demand system for the case study}
	\label{fig:flowsheet}
\end{figure*}

In iron \& steel manufacturing, the oxygen energy system is mainly made up of three subsystems (as illustrated in \figref{fig:flowsheet}): an oxygen generating system (OGS), an oxygen storage system (OSS) and an oxygen user system (OUS). OGS has a set of air separation units (ASUs) that simultaneously separates atmospheric air into pure gases, typically including nitrogen, oxygen and argon. The liquid oxygen pours into a tank, and the gaseous oxygen flows into the supply network connected with users. OSS includes the pipeline network that supplies high-pressure gas to each user and some oxygen gasholders units (GSUs) working as a buffer between supply and demand. OUS mainly consists of ironmaking and steelmaking workshops and other small users, in which they consume gaseous oxygen at different pressure. During the routine production, the ironmaking and steelmaking processes consume most of the oxygen, about 40\% for ironmaking and 50\% for steelmaking. In an iron \& steel plant, the decision support system named energy management system (EMS) focuses on determining optimal operational parameters and balancing the supply from OGS and the demand in OUS.

During the oxygen distribution horizon, we assume that: (1) Nitrogen, argon and other gases are sufficient. The distribution of oxygen cannot affect the EMS distributing other gaseous energy from OGS. (2) The OGS works under a normal condition and no machine breakdown happens. (3) The liquid oxygen is a type of recoverable energy and is only used when any ASU in OGS breakdown or the demand in OUS surges. Most of the efforts to handle oxygen distribution seek optimal operations and assume all parameters associated with manufacturing and energy are fully known. However, the balance between oxygen supply and demand often suffers from environmental and systematic uncertainty because of unforeseen events or unobservable factors, e.g. dynamic arrival of materials, duration variation of processes. In practice, it means that the energy distribution decisions made by the EMS may not fulfil the pre-designed balance. Therefore, how to distribute oxygen under uncertainty is a critical issue for the EMS achieving safety, stable and smart management.

The remainder of this paper is organized as follows: In Section \ref{S:2}, we review energy optimization models in iron \& steel plants and the robust optimization (RO) technique under uncertainty. The deterministic model for optimal oxygen distribution and its uncertain model using the two-stage robust optimization (TSRO) approach is investigated in Section \ref{S:3}. Section \ref{S:4} and \ref{S:5} are devoted to developing a time series model and a process scheduling model to describe the uncertainty of oxygen demands, respectively. In Section \ref{S:6}, computational studies are conducted to verify the effectiveness of the proposed approaches. The final section provides some conclusions and perspectives.

\section{Literature review}
\label{S:2}
The iron \& steel manufacturing system integrates with the energy flow and material flow, each of them restricts and promotes the other one \cite{ma2019energy}. The gaseous oxygen provides energy to manufacturing processes in OUS to make reaction or processing, and the operations of OGS and OSS are guided by the consumption demand of manufacturing processes \cite{sun2020materialb}. Therefore, how to make a good trade-off between the manufacturing system and the energy system is very important for modern iron \& steel plants to realize sustainable manufacturing \cite{o2014science}.  Many scholars and practitioners emphasized that this challenging issue can be tackled from the perspective of energy-efficiency and demand-side management \cite{merkert2015scheduling}. In the meantime, continued advances in mathematical modelling and information technologies also provide new opportunities to tackle the optimal distribution problem in the EMS of iron \& steel plants.

In recent studies, the energy distribution in an iron \& steel plant mainly focuses on making an optimal distribution with fixed demands by mixed-integer linear programming (MILP). Kong et al. \cite{kong2010milp}  and  Zhao et al. \cite{zhao2015milp,zhao2017optimal} proposed MILP models considering both operational and environmental factors to achieve simultaneously optimal by-product gases distribution in a Chinese iron \& steel plant. Recent studies also stressed the importance of optimal oxygen distribution in iron \& steel plants. Zhang et al. \cite{zhang2012model} formulated an oxygen optimal distribution problem with a linear programming model, in which the emission, pressure, liquefied, steam and setup costs are simultaneously minimized. Zhang et al.\cite{zhang2016milp} proposed a simultaneous multi-period model minimizing the oxygen emission which can make operational decisions involved on the supply side, e.g. oxygen generation rate of each ASU and compression rate of compressor.

Moreover, some studies also pointed out the optimal energy distribution can be achieved by exploring the flexibility on the demand side from the integrated perspective of energy and material flows \cite{gahm2016energy}. Zeng et al. \cite{zeng2018novel} proposed an improved MILP-based optimal distribution model by introducing decision variables named consuming rates of byproduct gases. Nolde and Morari \cite{nolde2010electrical} proposed a MILP-based scheduling solution to minimize the penalty caused by electricity over-and under consumption. Castro et al. \cite{castro2013resource} put forward a resource–task network to model an energy-constrained scheduling problem in iron \& steel plants, and studied the impact of fluctuating energy prices on operations and economic benefits. Hadera et al. \cite{hadera2015optimization} developed a  continuous-time MILP model simultaneously optimize a steelmaking scheduling and an electricity load commitment problem. However, it is interesting to note that the optimal oxygen distribution model needs to consider both continuous and discrete processes, each of which has different characteristics of demands. For an ironmaking process, its oxygen demand can be controlled by a continuous variable named flow rate \cite{zhou2018data}. For a steelmaking process, its oxygen demands always are adjusted by implements rescheduling. For example, xu et al. \cite{xu2018modeling} proposed an optimal steelmaking scheduling method to reduce the fluctuation of oxygen demands.

Due to the uncertainty in manufacturing and energy system, recent studies have reported data-based optimization methods to achieve more efficient decisions. Because the gasholder serves as a buffer during the distribution procedure, Zhang et al. \cite{zhang2011optimal} and Zhao et al. \cite{zhao2012effective} improved a least square support vector machine (LS-SVM) approach to forecast gasholder level with noises. Zhao et al. \cite{zhao2016data} proposed a data-based predictive optimization (DPO) method which forecasted the objective intervals with a Gaussian kernel-based learning model in the offline stage and implemented a time horizon rolling  MILP model in the online optimization stage. Some practitioners also tried to solve the optimal oxygen distribution with data-driven technique. Han et al. \cite{han2016two} applied a Granular-Computing (GrC)-based model to forecast the oxygen/nitrogen requirements and developed a MILP-based optimization model to distribute oxygen. Han et al. \cite{ han2017optimized} used an LS-SVM to fit the relationship between energy load and electricity cost of ASU and proposed a particle swarm optimization (PSO) algorithm to find optimal distribution. Although those methods can capture uncertainties, they only make optimal decisions in hindsight after uncertain factors are realized. Therefore, how to describe uncertainty in the realistic energy system (i.e. varied demands caused by production system) and integrate it with an optimal distribution model is a practical issue for making decisions in foresight. 

To tackle uncertainty in an optimal distribution model, we spontaneously employ stochastic programming (SP) or robust optimization (RO) to model its decision process. Comparing SO known probabilistic information, RO only assumes that the decision-maker has very little knowledge about the underlying uncertainty (except for its upper and lower bounds) and seeks an optimal solution that covers the worst-case cost within a well-synthesised uncertainty set. As early as 1973, Soyster \cite{soyster1973convex} began to study the robust linear programming problem. However, the RO technique had not been spread widely until  Ben-Tal and Nemirovski \cite{ben1998robust, ben1999robust} proposed a seminal theory. Recent literature \cite{bertsimas2011theory,gabrel2014recent,gorissen2015practical}  indicates RO has been widely used in transportation, supply chain, scheduling and has received growing attention as a modelling technique for energy system \cite{zugno2015robust,wang2015robust}. The essence of RO is to find out an optimal solution with strictly safe protection, but it is too conservative to apply in practice because all uncertain parameters are unlikely to reach the worst case at the same time. Therefore, this shortage motivates practitioners in the RO community on ways to seek alternative uncertainty sets with probabilistic guarantees, such as ellipsoidal, polyhedral set and budget-based set \cite{gorissen2015practical}. The budget-based uncertainty set proposed in \cite{bertsimas2004price} provides an intuitive interpretation of uncertainty for a risk-averse decision-maker who might flexibly adjust the level of risk aversion by tuning the so-called budget of uncertainty. It is important to note that this representation of uncertainty is general enough to account for interactions with energy management and can be easily interpreted by decision-makers.

Motivated by these practical industrial requirements and academic trends, we model oxygen distribution problem under uncertainty in the EMS with following main steps: 
\begin{itemize}
	\item Build an optimal oxygen distribution model considering flexible demands for continuous processes and discrete processes. The details will be stated in Section \ref{S:DO2D}.
	\item Derive a TSRO model with budget-based uncertainty set from the deeper insight into optimal oxygen distribution. The details will be stated in Section \ref{S:RO2D}.
\end{itemize}

Because oxygen users have different characteristics, we propose two methods to represent the uncertain oxygen demands of each distribution horizon:
\begin{itemize}
	\item Develop a time series model to forecast oxygen demands of continuous processes. The details will be stated in Section \ref{S:4}.
	\item Develop a process scheduling model to adjust oxygen demands of discrete processes. The details will be stated in Section \ref{S:5}.
\end{itemize}

\section{Optimal oxygen distribution model}
\label{S:3}
\subsection{Oxygen distribution in EMS}

To carry out an optimal decision for oxygen distribution, the EMS requires to balance supply and demand in each period ($\theta$) by setting operational parameters for both sides. On the supply side, the generation rate of each ASU ($\omega_{r,\theta}$) is controllable to produce oxygen volume. On the demand side, there are two approaches to adjusting oxygen demand. In the ironmaking process, oxygen-enriched hot air is continuously pumped into BF to produce molten iron \cite{nie2021numerical}.  Oxygen levels in the hot air can be adjusted according to oxygen supply. In the steelmaking process, multiple machines may work at the same time. Rescheduling simultaneous processing tasks is a beneficial way to avoid overloading or underloading \cite{xu2021reducing}. Therefore, we introduce two variables to define flexibility on the demand side.
\begin{itemize}
\item Adjustment rate ($\rho_{q}$), which represents the oxygen level in the hot air to the standard value.	
\item Scheduling scenario ($s$), which represents the schedule under a specific limitation of the multiple simultaneous processing tasks (so-called as capacity) in the steelmaking process.		
\end{itemize}

\begin{figure}[htbp!]
	\includegraphics[width=0.48\textwidth]{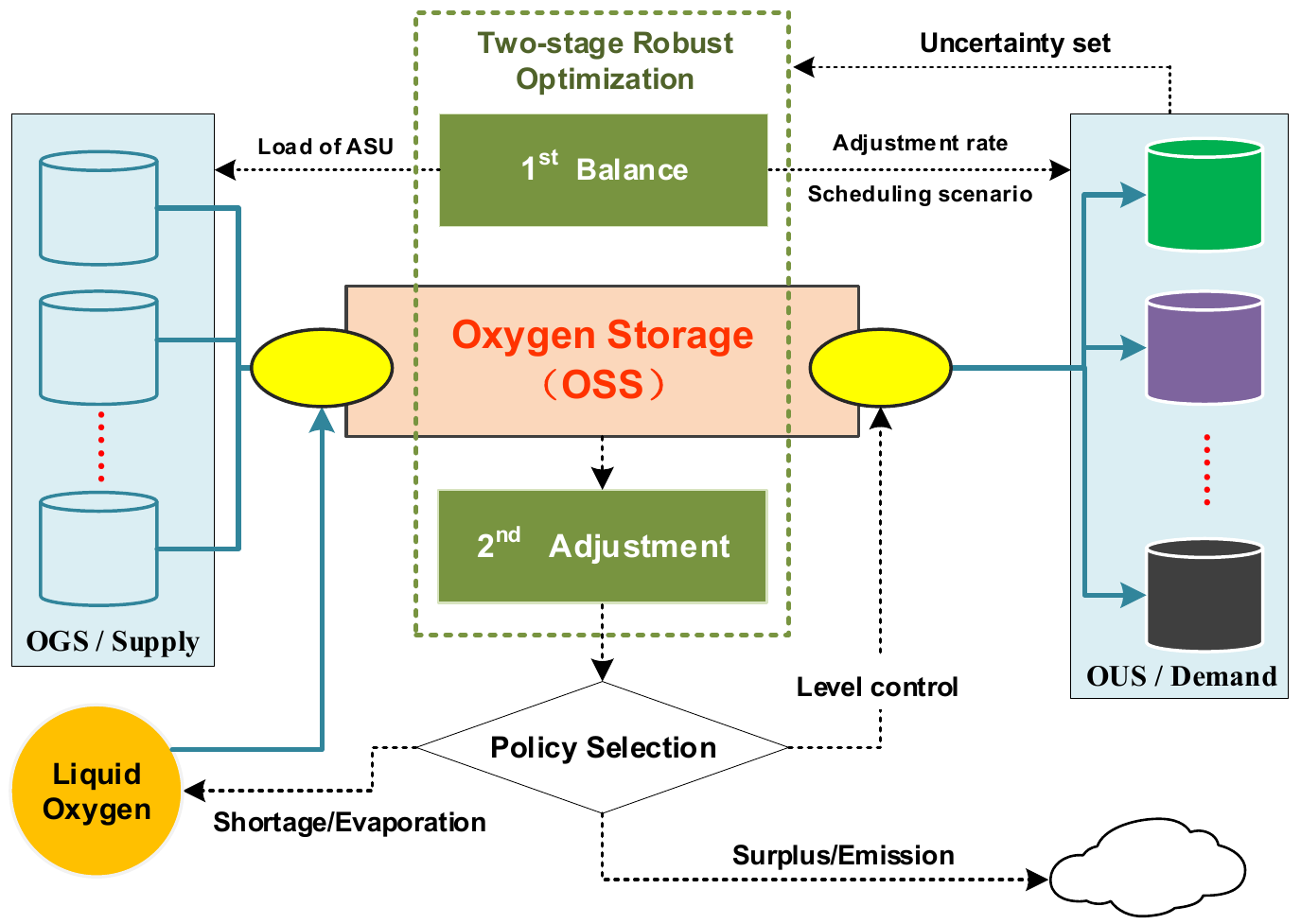}
	\caption{TSRO model for oxygen flexible distribution}
	\label{fig:policy}
\end{figure}

However, the oxygen demands in each period $d_{q,\theta}$ cannot be precisely estimated. When any demand in OUS is varied, the EMS needs to control the gasholder level to make a new balance between supply and demand. If these actions are ineffective, the gas in OSS will be in shortage or surplus and will cause extra cost. Considering the demands under uncertainty, we propose a TSRO model for oxygen flexible distribution as shown in \figref{fig:policy}. In the first stage, the operational parameters of the supply and demand sides are determined. In the second stage, an adjustment policy is chosen according to the status of the energy system.

In the following subsections, we first formulate a deterministic optimal oxygen distribution model involving flexible demands and develop a robust model which only knows limited information of demands.

\subsection{Deterministic model with flexible demands}
\label{S:DO2D}
\subsubsection{Objectives}

\begin{figure}[htbp!]
	\centering
	\includegraphics[width=0.45\textwidth]{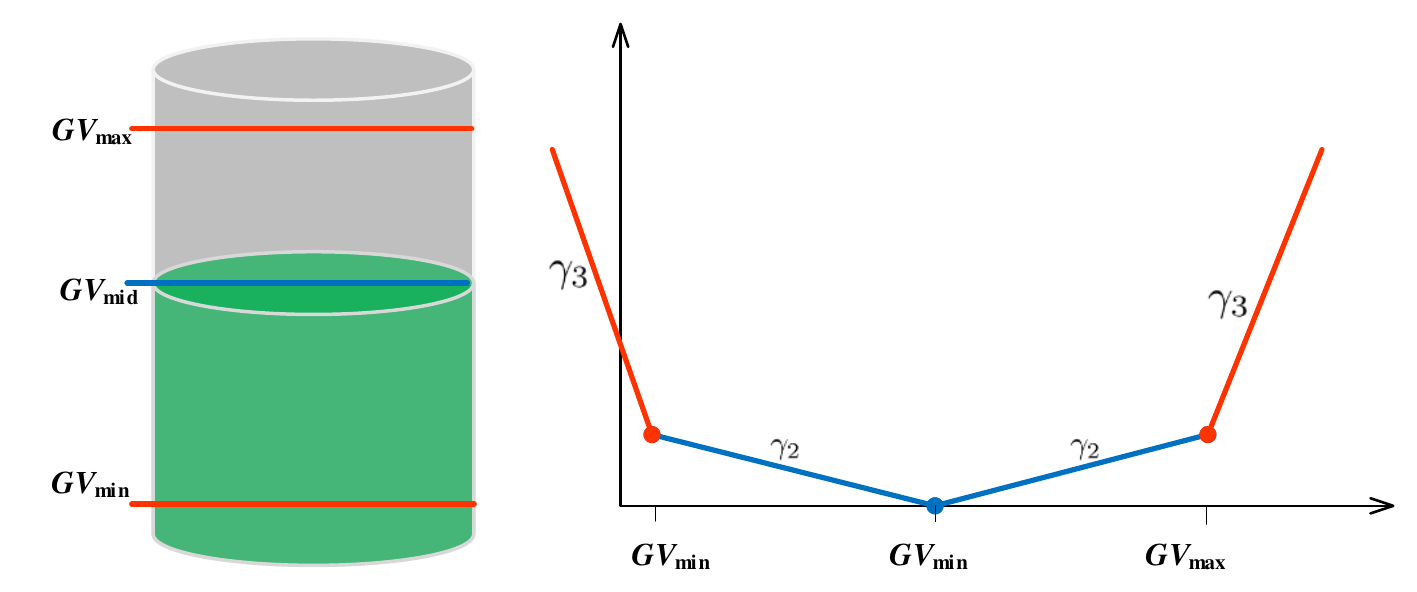}
	\caption{piecewise function for stability and penality}
	\label{fig:piece}
\end{figure}

In the balance stage, the EMS firstly needs to enlarge energy profit by maximizing the total load of OGS within the distribution horizon. In the adjustment stage, the gasholder level of OSS fluctuates frequently caused by varying demands. When the gasholder level is greater than $GV_{\max}$, the surplus will result in environmental pollution. When the gasholder level is less than $GV_{\min}$, the shortage will cause the extra cost by evaporating liquid oxygen to the gaseous one. Therefore, the middle level of the gasholder ($1/2$ between $GV_{\min}$ and $GV_{\max}$, as shown in \figref{fig:piece}) reflects the best comprise between shortage and surplus \cite{zhao2015milp}. Thus, the gasholder level should be retained approximately at the middle level to keep the best ability to reduce the risk of a gas system imbalance affects the manufacturing processes. In a word, the decision-maker needs to maximize the energy profit $f_1$ and minimize the gasholder deviation $f_2$ and the surplus/shortage penalty $f_3$:
\begin{equation}
	\begin{aligned}
		\max  f & = f_1\left(\omega\right) - f_2\left(\delta\right) - f_3\left(\epsilon\right) \\
		& = \gamma_1 \sum_{ r \in R} \sum_{\theta \in \Theta} \omega_{r,\theta}
		- \gamma_2\sum_{\theta \in \Theta}\left(\delta_{\theta}^{+} + \delta_{\theta}^{-} \right) \\
		& -\gamma_3 \sum_{\theta \in \Theta}\left( \epsilon_{\theta}^{+} +\epsilon_{\theta}^{-} \right) \\
	\end{aligned}
	\label{eq:oxy-obj}
\end{equation} 
where $\delta_{\theta}^{-}$ and $\delta_{\theta}^{+}$ respectively denote the deviation under and over $GV_{\rm{mid}}$, $\epsilon_{\theta}^{-}$ and $\epsilon_{\theta}^{+}$ respectively denote the volume of shortage and surplus.
\subsubsection{Constraints}
\begin{enumerate}[1)]
	\item Operational constraints: In the balance stage, the generation rate of ASU $r$  on the supply side and the adjustment rate of user $q$ on the continuous demand side are limited within the specified intervals (Eq.\ref{eq:Gminmax} , \ref{eq:Rminmax}). And only one scenario of the demands of discrete users can be selected (Eq. \ref{eq:scen}).
	\begin{equation}
		\omega_{r}^{\min} \leq \omega_{r,\theta} \leq \omega_{r}^{\max} \qquad \forall r \in R. \wedge \theta \in \Theta
		\label{eq:Gminmax}
	\end{equation} 
	\begin{equation}
		\rho_q^{\min} \leq \rho_{q} \leq \rho_q^{\max}   \qquad \forall q \in Q^N
		\label{eq:Rminmax}
	\end{equation} 	
	\begin{equation}
		\sum_{s \in S}y_s  =1     \qquad \forall y_s  \in \left\lbrace 0,1 \right\rbrace 
		\label{eq:scen}
	\end{equation} 
	
	\item  Ramp constraints: In the balance stage, the deviation of the generation rate of a ASU between two periods needs to be limited,
	$\left| \omega_{r,\theta} - \omega_{r,\theta-1} \right| \leq \hat{\omega}_{\max}
	$.
	Since $\left| \omega_{r,\theta} - \omega_{r,\theta-1}\right| $ is a non-linear function, we introduce two non-negative auxiliary variables ($\omega^{+}_{r,\theta}, \omega^{-}_{r,\theta}$), and let them satisfy:
	\begin{equation}
		\left\{
		\begin{array}{l}
			\omega_{r,\theta} - \omega_{r,\theta-1} = \omega^{+}_{r,\theta} - \omega^{-}_{r,\theta} \quad \forall  r \in R \wedge \theta \in \Theta\setminus\left\lbrace 1\right\rbrace   \\					
			\omega^{+}_{r,\theta} \geq 0, \, \omega^{-}_{r,\theta} \geq 0 \quad \forall  r \in R ,  \theta \in \Theta\setminus\left\lbrace 1\right\rbrace   \\			  
		\end{array} \right.
		\label{eq:Gaux}
	\end{equation}
	
	Let $\omega^{+}_{r,\theta} + \omega^{-}_{r,\theta} = \left| \omega_{r,\theta} - \omega_{r,\theta} \right|$, the ramp constraints can be reformulated as follows: 
	\begin{equation}
		\omega^{+}_{r,\theta} + \omega^{-}_{r,\theta} \leq  \hat{\omega}_{\max}  \qquad \forall  r \in R ,  \theta \in \Theta\setminus\left\lbrace 1\right\rbrace
		\label{eq:Gramp}
	\end{equation}	
	
	\item Flexible demands: the total energy demand in each period is determined by the adjustment rates of continuous users and the selected scenario of discrete users. 
	
	\begin{equation}
		d_{\theta}(\rho, z)  =  \sum_{q \in Q^N} \rho_{q} d_{q,\theta} + \sum_{q \in Q^E} y_{s} d_{q,\theta}^{s}   \quad \forall \theta \in \Theta		
		\label{eq:dmd}
	\end{equation} 
	
	\item Balance constraints: the oxygen generated at the end of the $(\theta-1)^{st}$ period flow into OSS by the beginning of the $\theta^{th}$  period, that is, all generated gas flows smoothly. According to the material equilibrium relationship proposed by Zhang et al. \cite{zhang2016milp}, the total of generated oxygen and user demand equals the storage change in the OSS, hence $GV_{\theta} - GV_{\theta-1} =  \sum_{r \in R} \omega_{r,\theta} - d_{\theta}(\rho, z)$. Because surplus demand causes emission and shortage of oxygen cause stock-out, we introduce two non-negative variables $\epsilon_{\theta}^{+}, \epsilon_{\theta}^{-} $ to denote the volume of surplus and shortage, respectively. Then, we can re-express  the gas balanced relationship with the following closed-form equation:	
	
	\begin{equation}
		\begin{aligned}
		GV_{0} + \sum_{r\in R}\sum_{\tau \in \Theta_{\theta}} \omega_{r,\tau} -& \sum_{\tau \in \Theta_{\theta}} \left( \epsilon_{\theta}^{+}-\epsilon_{\theta}^{-} \right) \\& =  \sum_{\tau \in \Theta_{\theta}} d_{\tau}(\rho, z)  
		 \quad \forall \theta \in \Theta
		\end{aligned}		
		\label{eq:Gdmd}
	\end{equation} 		
	
	\item Capacity constraints: commonly, the pressure of OSS is limited by its maximum and minimum value. According to the conversion relationship proposed by Zhang et al. \cite{zhang2016milp}, we define the safety volume of the OSS within a specific range, $\left[GV_{\min},  GV_{\max} \right] $. Thus, the capacity restrictions are represented in Eq.\eqref{eq:GVmin} and \eqref{eq:GVmax}. 	
	\begin{equation}
		\begin{aligned}
		GV_{0} +  \sum_{r\in R}&\sum_{\tau \in \Theta_{\theta}}   \omega_{r,\tau}   - \sum_{\tau \in \Theta_{\theta}} d_{\tau}(\rho, z)  \\ 
		 & -  \sum_{\tau \in \Theta_{\theta}} \left(\epsilon_{\theta}^{+}-\epsilon_{\theta}^{-} \right) 
		 \geq  GV_{\min}
		\forall  \theta \in \Theta	
		\end{aligned}	
		\label{eq:GVmin}
	\end{equation}  	 
	\begin{equation}
		\begin{aligned}
		GV_{0} + \sum_{r\in R}&\sum_{\tau \in \Theta_{\theta}} \omega_{r,\tau} - \sum_{\tau \in \Theta_{\theta}} d_{\tau}(\rho, z)  \\
		 & - \sum_{\tau \in \Theta_{\theta}} \left( \epsilon_{\theta}^{+} -\epsilon_{\theta}^{-} \right)
		\leq GV_{\max}
		\quad \forall \theta 	\in \Theta	
		\end{aligned}	
		\label{eq:GVmax}
	\end{equation} 
	
	\begin{equation}
		\epsilon_{\theta}^{+} \geq 0, \quad \epsilon_{\theta}^{-} \geq 0
		\quad \forall \theta 	\in \Theta		
		\label{eq:GVminax}
	\end{equation}

	\item Deviation definition: in the $\theta^{th}$ period, the deviation between the current gasholder level and the middle level can be expressed: $\left| GV_{\theta}-GV_{\rm{mid}} \right| $. Because the deviation is minimized in the objective function \eqref{eq:oxy-obj}, the deviation under and over the middle level can defined with following inequalities, respectively.  
	\begin{equation}
		\begin{aligned}
			\delta_{\theta}^{+} \geq GV_{0} + \sum_{r\in R}\sum_{\tau \in \Theta_{\theta}}& \omega_{r,\tau}   - \sum_{\tau \in \Theta_{\theta}} \left( \epsilon_{\theta}^{+}-\epsilon_{\theta}^{-} \right) \\ &- \sum_{\tau \in \Theta_{\theta}} d_{\tau}(\rho, z)-  GV_{\rm{mid}}
		\end{aligned}
		\label{eq:delta1}
	\end{equation} 	
	\begin{equation}
		\begin{aligned}
			\delta_{\theta}^{-} \geq  GV_{\rm{mid}} -\sum_{r\in R}\sum_{\tau \in \Theta_{\theta}}& \omega_{r,\tau}  + \sum_{\tau \in \Theta_{\theta}} \left( \epsilon_{\theta}^{+}-\epsilon_{\theta}^{-}  \right)
			\\ &
			 + \sum_{\tau \in \Theta_{\theta}} d_{\tau}(\rho, z) - GV_{0}
		\end{aligned}
		\label{eq:delta2}
	\end{equation} 	
	\begin{equation}
		\delta_{\theta}^{+} \geq 0, \quad \delta_{\theta}^{-} \geq 0
		\quad \forall \theta 	\in \Theta		
		\label{eq:DVminax}
	\end{equation} 	
\end{enumerate}

\subsubsection{Abstract model}

According the objective function and the constraint equations, the optimal oxygen distribution model with flexible demands and unvaried parameters is formulated as Eq.(\ref{eq:oxy-model}).

\begin{equation}
	\begin{aligned}
		&\max \quad f= f_1\left(\omega\right) - f_2\left(\delta\right) - f_3\left(\epsilon\right)\\	
		&s.t. \qquad  (\ref{eq:Gminmax})-(\ref{eq:DVminax})
	\end{aligned}
	\label{eq:oxy-model}
\end{equation} 

For notational brevity, we use set $\mathcal{D}$ to replace constraints (\ref{eq:Gminmax})-(\ref{eq:Gramp},\ref{eq:GVminax},\ref{eq:DVminax}), we complete the abstract optimal oxygen  distribution (AO2D) model with dimensionality reduction. 
\begin{equation}
	\rm \left( AO2D \right) \left\{
	\begin{aligned}
		\max \quad & 
		f_1 \left( {\omega} \right) - f_2 \left({\delta} \right) -  f_3\left( \epsilon \right)   \\		
		s.t. \quad & \sum_{\tau=1}^{\theta} \left(   {\omega}_{\tau} + {d}_{\tau} \left(\rho, z\right)  +   u_{\tau}(\delta,\epsilon)  \right)   \leq {v_i}  \quad 
		\\&\forall \theta \in \Theta, i \in \mathcal{I}, 
		\omega \in \mathcal{D}^{\omega}, \\
		& \qquad \qquad (\rho,y) \in \mathcal{D}^{\rho\times z},  (\delta,\epsilon) \in \mathcal{D}^{\delta\times \epsilon}  \\	  
	\end{aligned} \right.
	\label{eq:AO2D}
\end{equation}
where $\omega_{\theta}$ respectively denote two-dimensional variables $\omega_{r,\theta}$, and ${d}\left(\rho, z\right)$ and $ u(\delta,\epsilon)$ respectively represents demand function and balance control function and $ \mathcal{I}$ denotes associated constraint sets. 

\subsection{Robust model with limited information}
\label{S:RO2D}
Observing the formulation of the simplified model AO2D (\ref{eq:AO2D}), it is a variant of the multi-period inventory management problem and also is a two or multi-stage decision process. In this paper, we assume the oxygen demand of each user is uncertain and closely follow the idea of RO proposed by Bertsimas and Sim \cite{bertsimas2004price}, Bertsimas and Thiele \cite{bertsimas2006robust}, and develop a TSRO model involving the worst case. Therefore, the dynamic decision process is divided into two stages: (1) in the first stage, we have to determine here-and-now variables (i.e. generation rate $\omega_{r,\theta}$, adjustment rate $\rho_{q}$ and scenario selection $y_{s}$) in foresight; (2) in the second stage, wait-and-see variables are determined in hindsight where some uncertain demands are revealed. To make the worst-case RO be less-conservative, we also introduce the budget-based uncertainty set that minimizes the deviation to the best profit in a certain scenario. Taking these cluing, we develop the TSRO model by following three steps:

\subsubsection{Step 1: construct uncertainty set of demand} 

In this study, we assume that oxygen demand $d_{\theta}$ is a uncertain variable of which only the nominal value and bound are known. Let $[\bar{d}_{\theta} - \hat{d}_{\theta}, \bar{d}_{\theta} + \hat{d}_{\theta}]$, where $\bar{d}_{\theta}$ is the nominal demand and $\hat{d}_{\theta}$ is the maximal deviation.Then, we defined a scaled factor of the uncertain variable, $\xi_{\theta} = (\hat{d}_{\theta}  -\bar{d}_{\theta} ) /\hat{d}_{\theta}$, which falls into the closed range $[-1,1]$. Therefore, given a nominal value and a deviation, we can convert the uncertain variable into the function varied by a scaled factor, and impose a budget-based uncertainty set into the constraints within each time period $\theta$ as follows:

$$ \mathcal{U} = \left\lbrace \xi :  |\xi_{\theta}| \leq 1, \sum_{\tau=1}^{\theta} \xi_{\theta} \leq \Gamma_{\theta}, \forall  \theta \in \Theta  \right\rbrace  $$

With the predefined budget, we can avoid large deviations caused by many cumulative demands, obtain a “reasonable worst-case” solution. Its main assumption is that the  deviation increases with the number of periods $\theta$ considered, but the cumulative value should not exceed the budget  $\Gamma_{\theta}$. Thus, it is possible to control the degree of the conservativeness of the robust solution by varying the budget of uncertainty, obtaining a reasonable trade-off between optimality and robustness.

\subsubsection{Step 2: formulate robust counterpart model} 

According to the definition of RO theory, we formulate the RC model with a specified uncertainty set. Firstly, to determine the RC for each related inequalities, we replace
$d_{\theta}$ by $d_{\theta}= \bar{d}_{\theta} + \xi_{\theta}\hat{d}_{\theta} $ in model AO2D (\ref{eq:AO2D}), and guarantee that any solution is “immunized” against all realization of the uncertain demand must hold true for all $\xi$ over the budget-uncertainty set. Thus, we reformulate inequality constraints in model AO2D as follows:
\begin{equation}
	\begin{aligned}
		GV_0 + \sum_{\tau=1}^{\theta}\left( {\omega}_{\tau} - \left[ \bar{d}_{\tau} \left(  \rho, z\right) + \xi_{\theta} \hat{d}_{\tau} \left(  \rho, z\right) \right]  -  u_{\tau}(\delta,\epsilon)  \right)  \leq {v_i}  \\  \forall \theta \in \Theta , i \in \mathcal{I} , \xi \in \mathcal{U}
		\label{eq:GdmdX}
	\end{aligned}
\end{equation} 	

Based on the deterministic distribution model, we firstly formulated  a TSRO model with the worst case as follows:
\begin{equation}
	\rm \left( TSRO \right) \left\{
	\begin{aligned}
		&{\max_{\omega,\rho, z}} \quad   f_1(\omega) - {\min_{\xi \in \mathcal{U} }} \left\lbrace \Phi\left( \delta,\epsilon, \xi  \right) \right\rbrace  	 \\
		&\Phi\left(\xi \right) = \underset{\delta, \epsilon}{\min} \quad f_2(\delta) + f_3(\epsilon)   \\
		& s.t. \sum_{\tau=1}^{\theta}\left(u_{\tau}\left(\delta,\epsilon\right) + \hat{d}_{\tau}\left(\rho, z \right)\xi_{\tau}  \right)   \geq  GV_0  
		\\& \qquad\qquad \qquad +\sum_{\tau=1}^{\theta}\left({\omega}_{\tau} -\bar{d}_{\tau} \left(  \rho, z\right)   \right) -{v_i} \\
		&\qquad\qquad\qquad \qquad  \forall \theta \in \Theta, i \in \mathcal{I}, \xi \in \mathcal{U} 
	\end{aligned} \right.
	\label{eq:TSROWC}
\end{equation}

\subsubsection{Step 3: reformulate tractable model } 

The above formulation of TSRO  be recast via duality , which is suitable for directly using off-the-shelf optimizer, e.g., CPLEX, GUROBI, SCIP. For the sake of brevity, we show the reformulation procedure for the uncertainty constraint (\ref{eq:GdmdX}) according to its worst-case reformulation, which corresponds to maximize its left-hand side over the uncertainty
set $\mathcal{U}$. For the $\left( i,\theta\right)^{th}$ pair of uncertain constraints, the worst-case reformulation is as follows:

\begin{equation}
	\sum_{\tau=1}^{\theta}\left(  \omega_{\tau}  - \bar{d}_{\tau} ( \rho, z)- u_{\tau}(\delta,\epsilon) \right)   + 
	\Delta_{\theta}(\xi)
	\leq {v_i}
	\label{eq:GdmdY}
\end{equation}


Given a fixed vector of the pair $(\rho,z)$, the inequality (\ref{eq:GdmdY}) amounts to solving the two auxiliary linear programming problems: 
\begin{equation}
	\Delta_{\theta}(\xi) = \left\{
	\begin{aligned}
		&{\max} \quad  \sum_{\tau=1}^{\theta}\hat{d}_{\tau}(\rho, z) \xi_{\tau} \\
		& s.t. \quad  \sum_{\tau}^{\theta} \xi_{\tau} \leqslant \Gamma_{\theta}  \quad \forall \theta\\
		& \quad\quad -1\leqslant \xi_{\theta} \leqslant 1,   \quad \forall \theta
	\end{aligned} \right.
	\label{eq:AuLP}
\end{equation}  

Because the auxiliary linear programming problem \eqref{eq:AuLP} is feasible and bounded, by strong duality the optimal objective of this problem is equal to the optimal objective of its duality.  Reinjecting the dual of the auxiliary problem ($\ref{eq:AuLP}$) in models TSRO \eqref{eq:TSROWC}, we obtain the following tractable optimization problem \eqref{eq:TSRO}, where constraint (\rm{\romannumeral1}) imposes a set protective variables, (\rm{\romannumeral2}-\rm{\romannumeral3}) restrict the dual variables.

\begin{figure*}
\begin{equation}
	\rm \left( TSRO \right) \left\{
	\begin{aligned}
		&{\min_{\omega,\rho,y,\delta,\epsilon \in \mathcal{D} }} \quad   f_1 \left( {\omega} \right) - f_2 \left({\delta} \right) -  f_3\left( \epsilon \right)  & \\
		& s.t.\, \sum_{\tau=1}^{\theta}\left(\omega_{\tau}  - \bar{d}_{\tau} ( \rho, z)- u_{\tau}(\delta,\epsilon) \right)   + \beta_{\theta}\Gamma_{\theta} + \sum_{\tau=1}^{\theta} \alpha_{\theta,\tau} \leq {v_i}  
		\quad \forall \tau,\theta & (\rm{\romannumeral1})\\	
		& \qquad\quad \beta_{\theta} + \alpha_{\tau,\theta} \geq \hat{d}_{\tau} ( \rho, z) \qquad \forall \theta, \tau \leq \theta & (\rm{\romannumeral2})\\	
		& \qquad\quad \alpha_{\tau,\theta} \geq 0, \beta_{\theta} \geq 0   \qquad \forall \theta, \tau \leq \theta  & (\rm{\romannumeral3})
	\end{aligned} \right.
	\label{eq:TSRO}
\end{equation} 
\end{figure*}

\begin{figure}[htbp!]
	\centering
	\includegraphics[width=0.25\textwidth]{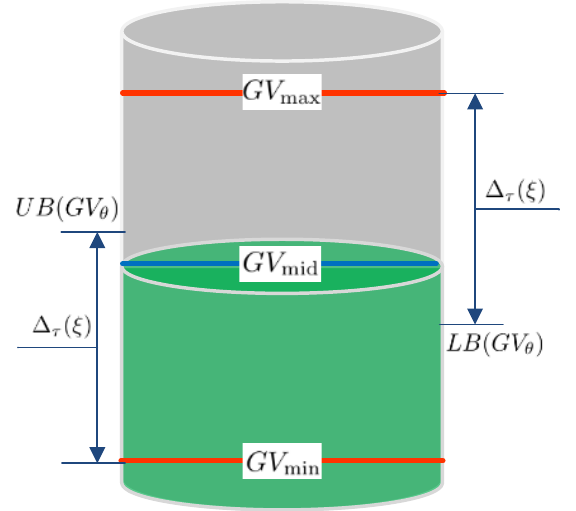}
	\caption{Upper and lower bound of gasholder level}
	\label{fig:bound}
\end{figure}

From Eq.\eqref{eq:GVmin},\eqref{eq:GVmax} and \eqref{eq:GdmdY}, we can deduce the upper and lower bound of gasholder level,
${\rm UB}\left( GV_{\theta} \right) = GV_{\max} - \Delta_{\theta}(\xi),  {\rm LB}\left( GV_{\theta} \right) = GV_{\min} + \Delta_{\theta}(\xi)$.
If ${\rm UB}\left(GV_{\theta}\right)<{\rm LB}\left(GV_{\theta}\right) $, TSRO \eqref{eq:TSRO} is supposed to be infeasible, as shown in \figref{fig:bound}. When ${\rm UB}\left(GV_{\theta}\right) \geq {\rm LB}\left(GV_{\theta}\right)$, $\Delta_{\theta}(\xi) \leq\left(GV_{\max}-GV_{\min}\right)/2$. Therefore, we need to bound $\Gamma_{\theta}$ when implement TSRO in practice.

\section{Interval demands forecasting for continuous processes}
\label{S:4}

Because the oxygen in the continuous processes consumed smoothly, we assume the nominal demand and its maximum deviation in each period are equal. That is, $\bar{d}_{q,\theta}= \bar{d}_{q}/\left| \Theta \right| $ and  $\hat{d}_{q,\theta}= \hat{d}_{q}/\left| \Theta \right| $, where $\bar{d}_{q}$ and $\hat{d}_{q}$ are the nominal demand and the maximum deviation user $q$ over the distribution horizon. In this section, we develop a time series model to forecast the oxygen demand of the next distribution horizon.

\subsection{Time-series forecasting}
In iron \& steel manufacturing plants, energy data is collected as time series which is an observation of some underlying dynamical process of interest. Some practitioners applied these day-marked energy consumption data and perform a model-based analysis to forecast future values based on previous observations \cite{sun2013plant,pena2019optimal}. 


Commonly, time series tends to exhibit high correlations induced by the temporal structure in the data. Consider a time series of oxygen demands $d_t$, $t\in T$ parametrized by the set $T$,  the goal of the time-series forecasting model is to find a function $f(\cdot)$ based on a state-space model ($y_{t}$) \cite{girard2003gaussian}. 
\begin{equation}
	d_{t}=f\left(y_{t}\right) +  \varepsilon_{t} = f\left(d_{t-1},\dots, d_{t-TL} \right) +  \varepsilon_{t} 
	\label{eq:ts}
\end{equation} 
where $TL$ is represents the size of time lags, $t$ indicates the moment in the past we focus on and $\varepsilon_{t}$ is a white noise. Then, we assume already have obtained a training input $\bm{\mathbf{y}} = (y_{t-1},\dots, y_{t-N})$ with size of $N$. 

Generally speaking, the time series of energy users has different patterns because they have different process characteristics which mainly include: seasonality, trend, noise, change point. Due to limited assumptions over the underlying function, traditional linear time series models based on a mathematical formula and techniques are hard to model complex pattern.GP (see \cite{rasmussen2003gaussian} for a thorough background) permits a convenient and natural framework for modelling a time series forecasting function $f(\cdot)$ because of its flexibility with prior knowledge and uncertainty quantification with estimated probabilistic interval.  


\subsection{GP-based time-series model}

GP is a stochastic process for $f(y_t)$  where for selected time points $t_1, t_2, \dots, t_{N}$ of which the probability distribution $\Pr\left[f(y_{t_1}),\dots, f(y_{t_N})\right] $ is multivariate Gaussian. To formulate a GPR model, we first define a mean function $\mu(y_t)=\mathbb{E}\left[f(y_t) \right]  $ and a covariance function (also referred to a kernel function) 
$k\left(y_{t_1}, y_{t_2}\right) ={\rm Cov}\left[f(y_{t_1}), f(y_{t_2}) \right] = \mathbb{E}\left[(f(y_{t_1})-\mu(y_{t_1}))(f(y_{t_2})-\mu(y_{t_2}))\right]$,and denote the GP by

\begin{equation}
	f(y_t) \sim \mathcal{GP} \left(\mu(y_{t}), k\left(y_{t}, y_{t'}\right)  \right)
	\label{eq:gp}
\end{equation} 

For training set $\bm{\mathbf{y}} = \left\lbrace y_{t_1}, y_{t_2},\cdots, y_{t_N} \right\rbrace $, we can define following covariance matrices:

\begin{equation}
	\bm{K}\left(\bm{\mathbf{y}}, \bm{\mathbf{y}}\right) =  
	\begin{bmatrix}
		k(y_{t-1},y_{t-1}) & \cdots   & k(y_{t-1},y_{t-N}) \\ 
		\vdots & \ddots   & \vdots\\ 
		k(y_{t-N},y_{t-1}) &  \cdots  &  k(y_{t-N},y_{t-N})
	\end{bmatrix}
	\label{eq:cov1}
\end{equation} 

\begin{equation}
	\bm{{K}}\left(y_t, \bm{\mathbf{y}}\right) = \left[ k(y_{t},y_{t-1}),\cdots,k(y_{t},y_{t-N}) \right]
	\label{eq:cov2}
\end{equation} 

According to the Bayes' Theorem, we can infer the posterior distribution by defining the new joint distribution:

\begin{equation}
	\Pr\left(\begin{bmatrix}
		\bm{f(\bm{\mathbf{y}})} \\
		f(y_t) 
	\end{bmatrix} \right) = \mathcal{N} \left( \begin{bmatrix}
		\bm{\mu(\mathbf{y})} \\
		\mu(y_{t}) 
	\end{bmatrix} , \begin{bmatrix}
		\bm{{K}}\left( \bm{\mathbf{y}}, \bm{\mathbf{y}} \right) &\bm{{K}} \left( \bm{\mathbf{y}}, y_t \right)  \\
		\bm{{K}}\left(y_t, \bm{\mathbf{y}} \right) & k\left(y_t, y_t\right) 
	\end{bmatrix}  \right)  
	\label{eq:joint}
\end{equation} 

Through some manipulations of the joint Gaussian distribution, the resulting posterior distribution over $f(y_t)$  is also a GP with mean and covariance
functions given by:
\begin{equation}
	\mu_{*}\left(y_t\right) = \bm{K}\left(y_t, \bm{\mathbf{y}} \right) \bm{K}\left(\bm{\mathbf{y}}, \bm{\mathbf{y}} \right)^{-1}\left( {f(\bm{\mathbf{y}})- \mu(\bm{\mathbf{y}})} \right) 
\label{eq:postm}
\end{equation} 
\begin{equation}
	k_{*}\left(y_t, y_t\right) =  k\left(y_t, y_t\right)-  \bm{{K}}\left(y_t, \bm{\mathbf{y}} \right)\bm{{K}}\left( \bm{\mathbf{y}}, \bm{\mathbf{y}} \right)^{-1}\bm{{K}}\left( \bm{{y}},y_t \right) 
\label{eq:postv}
\end{equation} 
While $\mu_{*}$  can be used as a point estimate and $k_{*}$ can be used as a interval estimate, we have a full distribution for $d_t$, allowing us to quantify the associated uncertainty in our predictions. In this paper, we only focus on one-step ahead predictions as a first step and leave multiple-step ahead forecasts for future work (more details please \cite{girard2003gaussian}).

Commonly, the prior mean $\mu(\bm{\mathbf{y}})$  of a GP often is set to 0. 
The covariance function of a GP specifies the correlation between any pair of outputs. This can then be used to generate a covariance matrix over our set of observations and predictions. Fortunately, there exist a wide variety of functions that can serve in this purpose \cite{rasmussen2003gaussian}, which can then be combined and modified in a further multitude of ways. This gives us a great deal of
flexibility in our modelling of functions, with covariance functions available to model periodicity, delay, noise and long-term drifts and other phenomena.

\section{Multi-scenario demands for discrete processes}
\label{S:5}

Because the oxygen in the discrete processes consumed non-smoothly, we cannot directly calculate the nominal demand and its maximum deviation in each period. Thus, we need to calculate the interval demands according to the schedule of steelmaking processes. That is, $\bar{d}_{q,\theta}= \bar{d}_{q}*TP_{q,\theta}/TP_{q}$ and  $\hat{d}_{q,\theta}= \hat{d}_{q} *TP_{q,\theta}/TP_{q}$, where $TP_{q}$ and $TP_{q,\theta}$ represent the total processing time allocated to user $q$ and the one in period $\theta$,  $\bar{d}_{q}$ and  $\hat{d}_{q}$ are estimated with GP-based approach in our previous work \cite{jiang2019gaussian}. In this section, we propose a scheduling method for the EMS to generate multiple scheduling solutions via defining different capacity limits. 

\subsection{Scheduling for steelmaking production}
\begin{figure}[htbp!]
	\centering
	\includegraphics[width=0.45\textwidth]{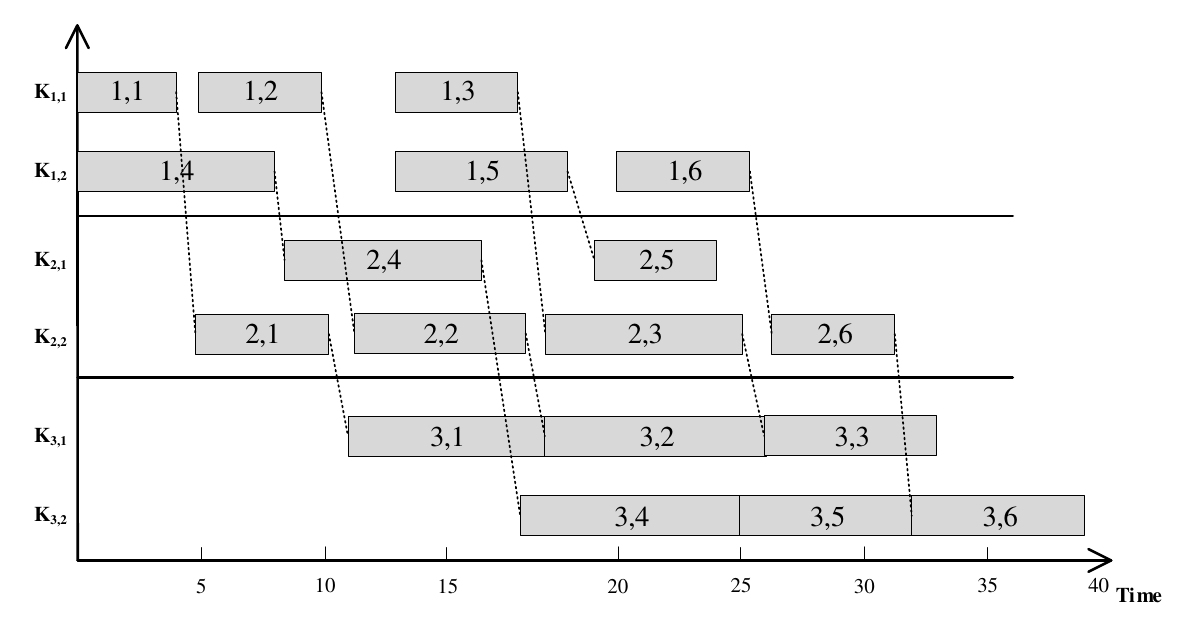}
	\caption{A Gantt chart for the feasible schedule executed in the steelmaking workshop}
	\label{fig:ganttdemo}
\end{figure}

\begin{figure}[htbp!]
	\centering
	\includegraphics[width=0.40\textwidth]{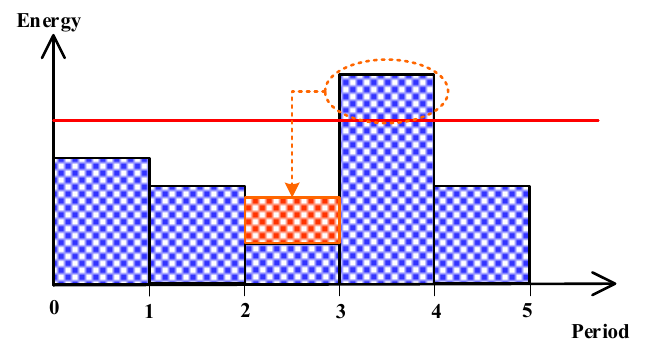}
	\caption{Peak clipping \& valley filling on energy demand shape}
	\label{fig:peakshift}
\end{figure}

Scheduling in steelmaking processes is linked with three types of stages: steelmaking,  refining and continuous casting, through which high-temperature molten iron from ironmaking processes is converted into solidified slabs. It focuses on sequencing a job set $J$ each of which has a set of tasks $A_j$, allocating an eligible machine for each task $A_{i,j}$, and optimizing some given criteria, such as make-span, waiting time, and earliness/tardiness. The task set of each job must follow the same order from steelmaking to  casting stage. In the casting stage, tasks must be performed in batch and without interruption because the liquid steel is solidified continuously within a period. Considering such procedural, temporal and resource constraints, the scheduling model of the steelmaking processes is always identified as a variant of the hybrid flow shop-scheduling problem \cite{tan2013two,xu2018modeling}. Then, a Gantt chart containing several continuous blocks at the last stage can represent a feasible or optimal schedule executed in all stages (as shown in \figref{fig:ganttdemo}). Given a fixed schedule, tasks in the steelmaking stage consume a large amount of oxygen and often cause a non-smooth demand shape because of the discrete processes in scheduling. To reduce the impact on the oxygen distribution system, the demand shape needs to be shifted by peak clipping and valley filling (as shown in \figref{fig:peakshift}). 


\subsection{Preprocessing procedure}

\begin{figure}[htbp!]
	\centering
	\includegraphics[width=0.45\textwidth]{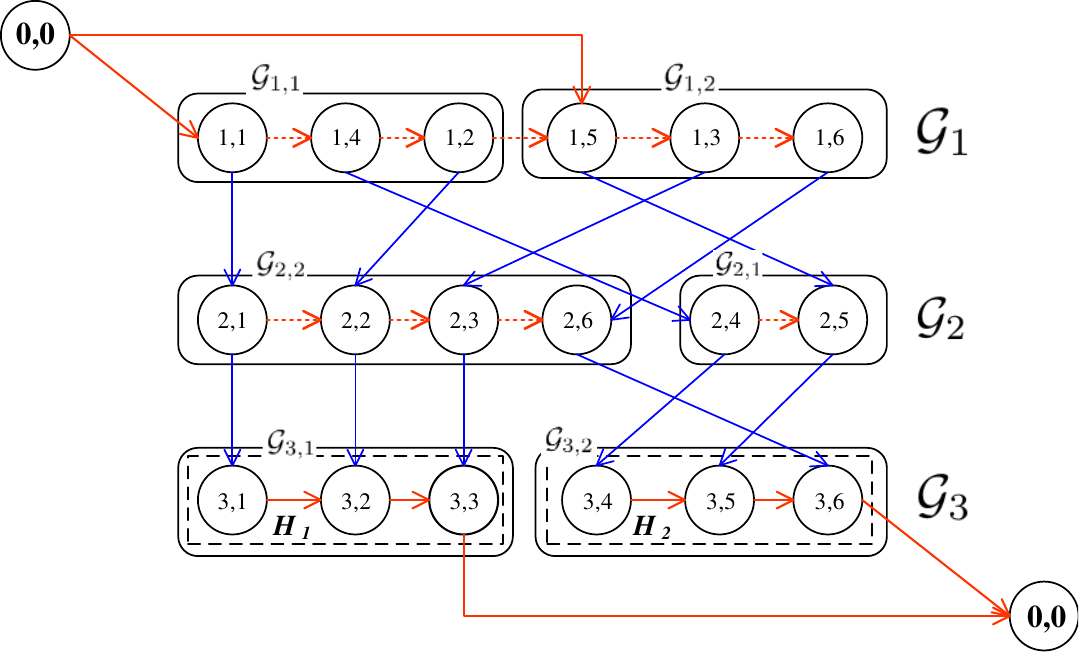}
	\caption{A directed graph $\mathcal{G}$ for a schedule released into a steelmaking workshop}
	\label{fig:digraph}
\end{figure}

Because the number of simultaneous processing tasks named capacity directly determines the oxygen consumption, we propose a capacity-constrained shifting method to adjust oxygen on the demand side. Depending on the problem size and its complexity, the scheduling methodologies involving resource constraints can be divided into two groups: monolithic and sequential approaches \cite{marchetti2009general}. The monolithic approach simultaneously solves the job assignment, sequencing problem and timing problem, while the sequential approach assumes that the job sequencing and machine allocation problems have already been solved. Considering the advantage on computational complexity, we adopt the sequential approach and represent the solution of the job sequencing and machine allocation with a directed graph $\mathcal{G}$ shown in \figref{fig:digraph}, where $A_{0,0}$ is a dummy task.

Since timing optimization problem is highly sensitive to the time horizon, we perform a preprocessing procedure on the directed graph and compute the earliest starting times ($ES_{g,j}$) and the latest finishing times ($LF_{g,j}$) for each task $A_{g,j}$. The time-window preprocessing benefits: 
\begin{itemize}
	\item Reducing the number of the time-indexed decision variables,
	\item Reducing the total number of constraints needed.
\end{itemize}
Given a fixed $\mathcal{G}$, we let:
\begin{itemize}
	\item $\mathcal{G}_g$ and $\mathcal{G}_{g,m}$ respectively represent tasks allocated in stage $g$ and the $m^{th}$ machine in $M_{g}$,
	\item $A_{gg,j}$ is the next task of $A_{g,j}$ and  $A_{g,jj}$ is the next task of $A_{g,j}$ on the same machine
\end{itemize}
Here, the two-pass method named forward-backward scheduling (stated in Algorithm \ref{alg:FBS}) is used to estimate $ES_{g,j}$ and $LF_{g,j}$.

\begin{algorithm}
	\caption{forward-backward scheduling}
	\label{alg:FBS} 
	\LinesNumbered 
	\KwIn{A directed graph $\mathcal{G}$ }
	\KwOut{Earliest starting Time($ES_{i,j}$) , Latest finish times($LF_{i,j}$)}%
	{/*Notes: $\bar{S}_{g,j}$ $\bar{C}_{g,j}$ represent estimated starting and finishing time of $A_{g,j}$ */} \\
	\text{/*Step 1. forward scheduling*/} \\
	Set $\bar{S}_{0,0}:=0$ \;
	\For{ $gg=1$ to $\left|G\right|$  }{
		
		\ForEach{$A_{gg,jj} \in \mathcal{G}_{gg}$}
		{
			$\bar{S}_{gg,jj}:=\bar{S}_{gg,j}+PT_{gg,j} + 1 $ \;
		}			
		\ForEach{$A_{gg,j} \in \mathcal{G}_{gg}$ }
		{
			$ES_{gg,j} := \max\left\lbrace \bar{S}_{g,j}+PT_{g,j}+TT_{g,gg} + 1, \bar{S}_{gg,j} 
			\right\rbrace$\;	
			$ \bar{S}_{gg,j} := ES_{gg,j} $\;	
		}
	}
	\text{/*Step 2. backward scheduling*/} \\
	$\bar{S}_{0,0} = T$ \\
	\For{ $g=\left|G\right|-1 $ to $1$}{
		
		\ForEach{$A_{g,j} \in \mathcal{G}_{g}$}
		{
			$\bar{C}_{g,j}=\bar{C}_{g,jj}-PT_{g,jj}-1$ \;
		}			
		\ForEach{$A_{g,j} \in \mathcal{G}_{g}$}
		{
			$LF_{g,j} := \max\left\lbrace \bar{C}_{gg,j}-PT_{gg,j}-TT_{g,gg}-1, \bar{C}_{g,j} 
			\right\rbrace$\;	
			$ \bar{C}_{g,j} := LF_{g,j}$\;	
		}
	}
\end{algorithm}

\subsection{Model formulation}
In this section, we developed a mixed integer linear programming (MILP) discrete-time models for the capacity-constrained steelmaking scheduling problem, which is based on the definition of binary variables ($x_{g,j_l}$) that describe the processing state of every task $A_{g_j}$ at each time point $l$.
\subsubsection{Objectives}
Considering productivity and energy loss, the objective of the steelmaking scheduling problem is to minimize the total of make-span ($C_{\max}$) and waiting times ($W_{\rm{tot}}$).
\begin{equation}
	\begin{aligned}
	\min F & =  C_{\max} + W_{\rm{tot}} \\ &= \max_{j\in J}{C_{\left| G \right| ,j}} +  \sum_{j \in J}\left( C_{\left| G \right| ,j} - C_{1,j} - \sum_{g=1}^{\left| G \right| -1 } PT_{g,j} \right) 
	\end{aligned}
	\label{eq:scc_obj}
\end{equation} 

\subsubsection{Constraints}
A feasible schedule executed in the steelmaking process need to satisfy following constraints.   
\begin{enumerate}[1)]	
\item Each task is be only started at exactly one time unit.
\begin{equation}
	\sum_{l=1}^{|L|-PT_{g,j}+1} x_{g,j,l} =1, \quad \forall  A_{g,j}
	\label{eq:time_unique}
\end{equation} 
\item The completion time of task $ A_{g,j} $ can be defined as:
\begin{equation}
	C_{g,j} = \sum_{l=1}^{|L|-PT_{g,j}+1} \left( l-1 \right) x_{g,j,l} + PT_{g,j}, \quad \forall A_{g,j} 
	\label{eq:scc-finish}
\end{equation} 
\item For every two connected tasks of the same job (i.e., $ A_{g,j}$ and $ A_{gg,j}$ ), the next operation cannot be started unless the previous one has been completed and delivered.
\begin{equation}
	C_{gg,j} - C_{g,j} \geq PT_{g,j}+TT_{g,gg} +1 , \quad \forall  A_{g,j},  A_{gg,j}
	\label{eq:scc_link}
\end{equation}
\item For every two connected task on the same machine (i.e., $A_{g,j}$ and $A_{g,jj}$ ), the next task cannot be started unless the previous one has been completed.
\begin{equation}
	C_{g,jj} - C_{g,j} \geq PT_{g,jj} +1, \quad \forall  A_{g,j}, A_{g,jj}
	\label{eq:scc_seq}
\end{equation}
\item The number of jobs simultaneously processed at time $t$ cannot be greater than the total number of available resource. $\mathcal{G}_O$ represents the stages consuming oxygen.
\begin{equation}
	\sum_{ A_{g,j} \in \mathcal{G}_O } x_{g,j,l}  \leq WL_{q}^{\max}, \quad \forall  l \in \left[ ES_{g,j},  LF_{g,j}\right]
	\label{eq:scc_capcity}
\end{equation} 
\item A minimum setup time at the casting stage must be guaranteed before the first job of a new batch starts.
\begin{equation}
	C_{\left| G \right|,(hh)^{'}} - C_{\left| G \right|,(h)^{''}} = PT_{\left| G \right|,(hh)^{'}} + SU_{hh}, \quad \forall h,hh
	\label{eq:scc_setup}
\end{equation} 
where $(\cdot)^{'}$ and $(\cdot)^{''}$ respectively denote the first and the last job in a batch.
\item The processing precedence constraints of jobs in the same batch in the last stage.
\begin{equation}
	C_{\left| G \right|,jj}-C_{\left| G \right|,j} = PT_{\left| G \right|,jj} + 1, \quad   \forall h\in H \wedge j,jj \in J(h)
	\label{eq:scc_continuity}
\end{equation} 
\end{enumerate}

\section{Computational studies}
\label{S:6}
In this section, we present computational studies to verify our proposed models and test their ability to effectively distribute oxygen under uncertainty. The TSRO and scheduling model were programmed by Pyomo \cite{hart2017pyomo}, a Python-based, open-source algebraic modelling language, and solved by "Gurobi 9.0" (\url{https://www.gurobi.com/products/gurobi-optimizer/}, the academic version with default settings), a state-of-the-art optimizer for mathematical programming. For the GP implementation, we adapted the open source scikit-learn package with version 0.24.1 \cite{pedregosa2011scikit}. All computational studies were executed on a personal computer with an Intel Core i7 processor (3.60 GHz), 16.0GB RAM and Windows 10 operating system.

Case studies addressed in this section were performed by applying our methods in an iron \& steel plant of China, in which there were two ASUs in OGS for supplying oxygen, five oxygen users in OUS including two ironmaking processes, two primary steelmaking processes and a small user from other processes. Computational studies were proceeded as following:
\begin{enumerate}[(1)]	
	\item Forecasted oxygen demands of each horizon with the GP-based time series, estimated the point and interval values for the TSRO model..
	\item Generated multi-scenario oxygen demands by rescheduling production tasks in the steelmaking processes.
	\item Analysed the proposed TSRO model for oxygen distribution under uncertainty including its budget levels, sensitivity and robustness.	  
\end{enumerate}

 In following studies, we collected statistic data of manufacturing and energy data from the EMS and randomly synthesized energy distribution optimization instances with different demand levels. All testing instances were shared on the Github webpage  (\url{https://github.com/janason/Energy/tree/master/O2}).

\subsection{Interval demands forecasting by GP}

To validate the forecasting model, we collected three months (90 days) of actual daily oxygen data. Three are three samples  each day, thus $90\times 3= 270$ samples in total. We used 235 samples for training and 35 samples for testing oxygen demand with one horizon ahead. The GP-based time series model was implemented with Scikit-Learn's class \textit{GaussianProcessRegressor}. We specified the prior mean to be \textit{zero} and prior covariance by passing a mixture \textit{kernel} object with a squared exponential and a white noise covariance function. We ran the forecasting model with ten random restarts, a maximum of 2,000 iterations and evaluated it with a metric named Mean Absolute Percentage Error (MAPE).

$$ MAPE=\frac{\sum_{t=1}^{N}f(y_t)-y_t}{N} \times 100 \% $$

\begin{figure}[htbp]
	\centering
	\subfigure[ironmaking process 1\# (MAPE = 1.2\%)]{
		\begin{minipage}[t]{1.0\linewidth}	
			\centering		
			\includegraphics[width=1.0\linewidth]{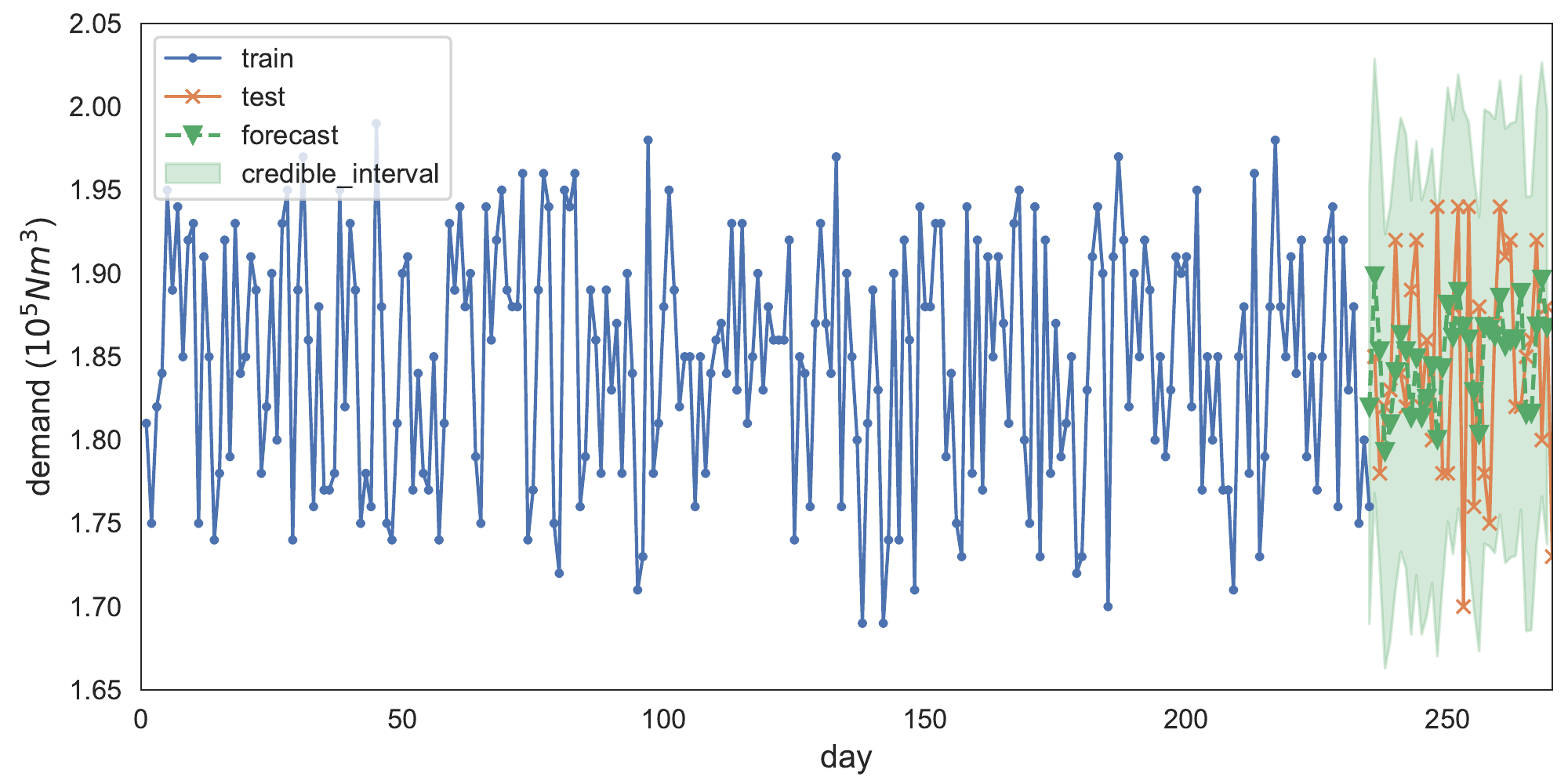}
		\end{minipage}%
	}\\
	\subfigure[ironmaking process 2\# (MAPE = 2.1\%)]{
		\begin{minipage}[t]{1.0\linewidth}	
			\centering		
			\includegraphics[width=1.0\linewidth]{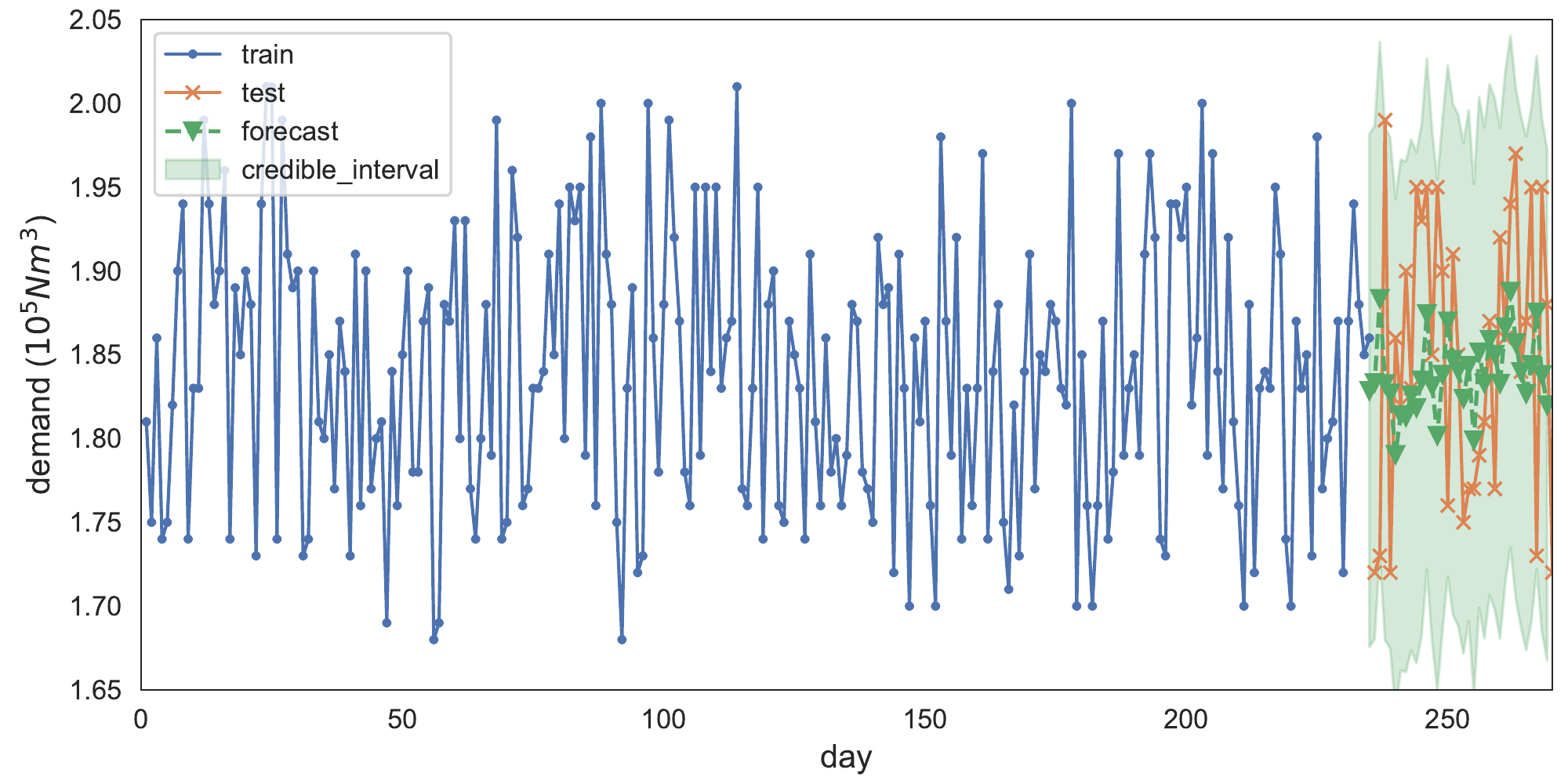}
		\end{minipage}%
	}%
	\caption{Probabilistic forecasts of oxygen interval demand}
	\label{fig:GPTS}
\end{figure}

To investigate the forecasting performance in a more intuitive way, we gave a graphical illustration of the probabilistic forecasting results for the two ironmaking processes. In \figref{fig:GPTS}, we respectively plot the observed and forecasted oxygen demands, represent their $95\%$ confidence intervals with shading boundaries.  We observe that almost all points in the testing set fall into the boundaries. Both MAPEs of two ironmaking processes implementing GP method are less than $2.5\%$, which means the average errors between actual and forecasted values is small and acceptable. After probabilistic forecasting of the next distribution horizon, the estimated point and interval information of oxygen demand can be obtained and provided for the TSRO model.

\subsection{Multi-scenario demands by capacity-constrained rescheduling}

As previously alluded to, there often exist multiple scenarios for steelmaking production scheduling, and each of them has different demands on oxygen. The iron \& steel plant has two steelmaking processes that consume oxygen in the meantime. In the first process, two machines are for dephosphorization (DP). In the second process, three machines are for decarbonization (DC). In this study, we selected three representative scheduling instances. Their multi-scenario energy demands were generated by setting two maximum capacities of the steelmaking processes: 
\begin{itemize}
	\item Scenario 1, the maximum capacity $WL_{1}^{\max} = 5 $, where the number of simultaneous processing tasks is not greater than 5.
	\item Scenario 2, the maximum capacity $WL_{2}^{\max} = 4 $, where the number of simultaneous processing tasks is not greater than 4.
\end{itemize}

\begin{figure}[htbp]
	\centering
	\subfigure[Scenario 1: $C_{\max}=704,W_{\rm{tot}}=304$]{
		\begin{minipage}[t]{1.0\linewidth}	
			\centering		
			\includegraphics[width=1.0\linewidth]{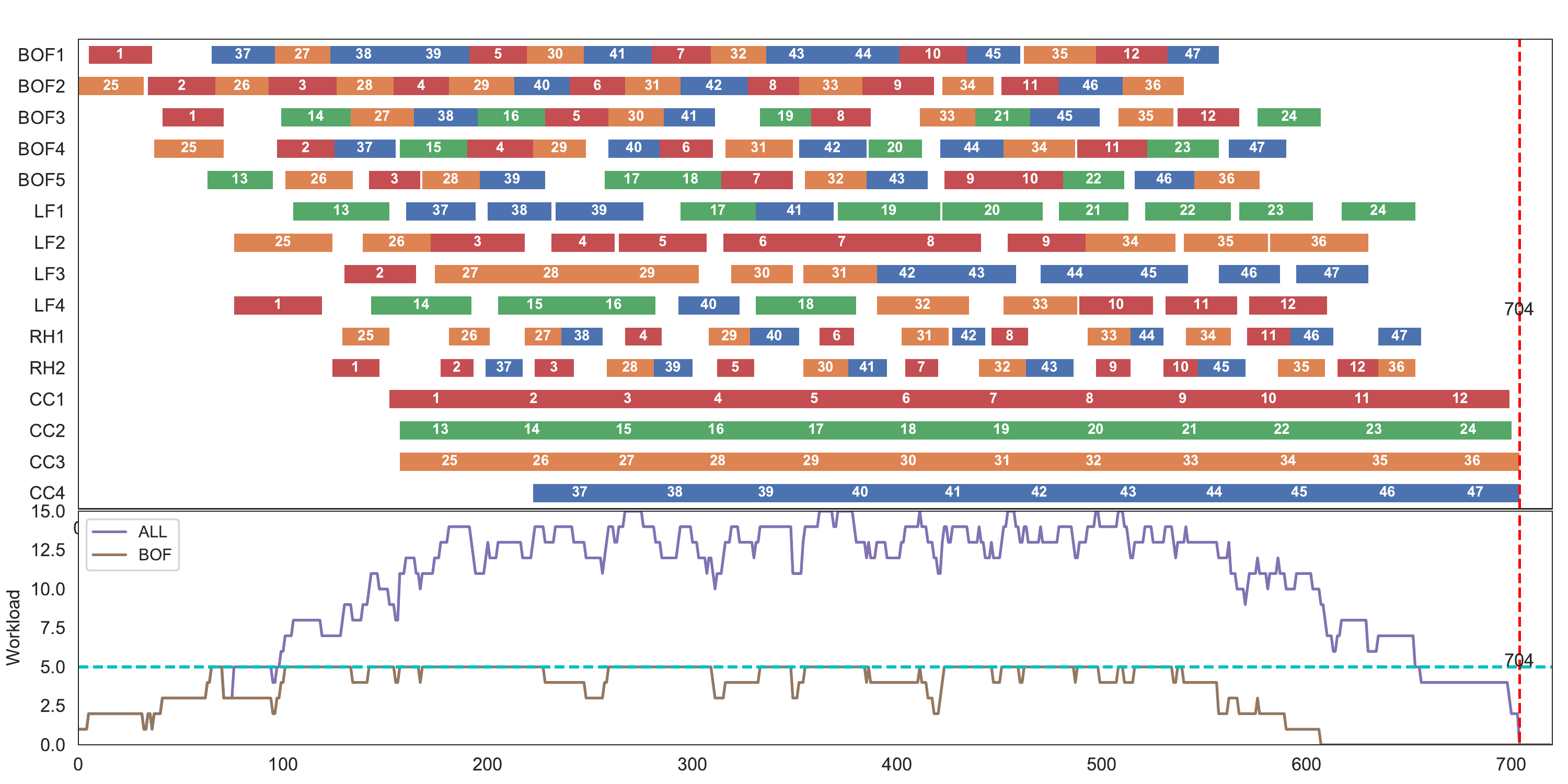}
		\end{minipage}%
	}\\
	\subfigure[Scenario 2:$C_{\max}=819,W_{\rm{tot}}=3642$]{
		\begin{minipage}[t]{1.0\linewidth}	
			\centering		
			\includegraphics[width=1.0\linewidth]{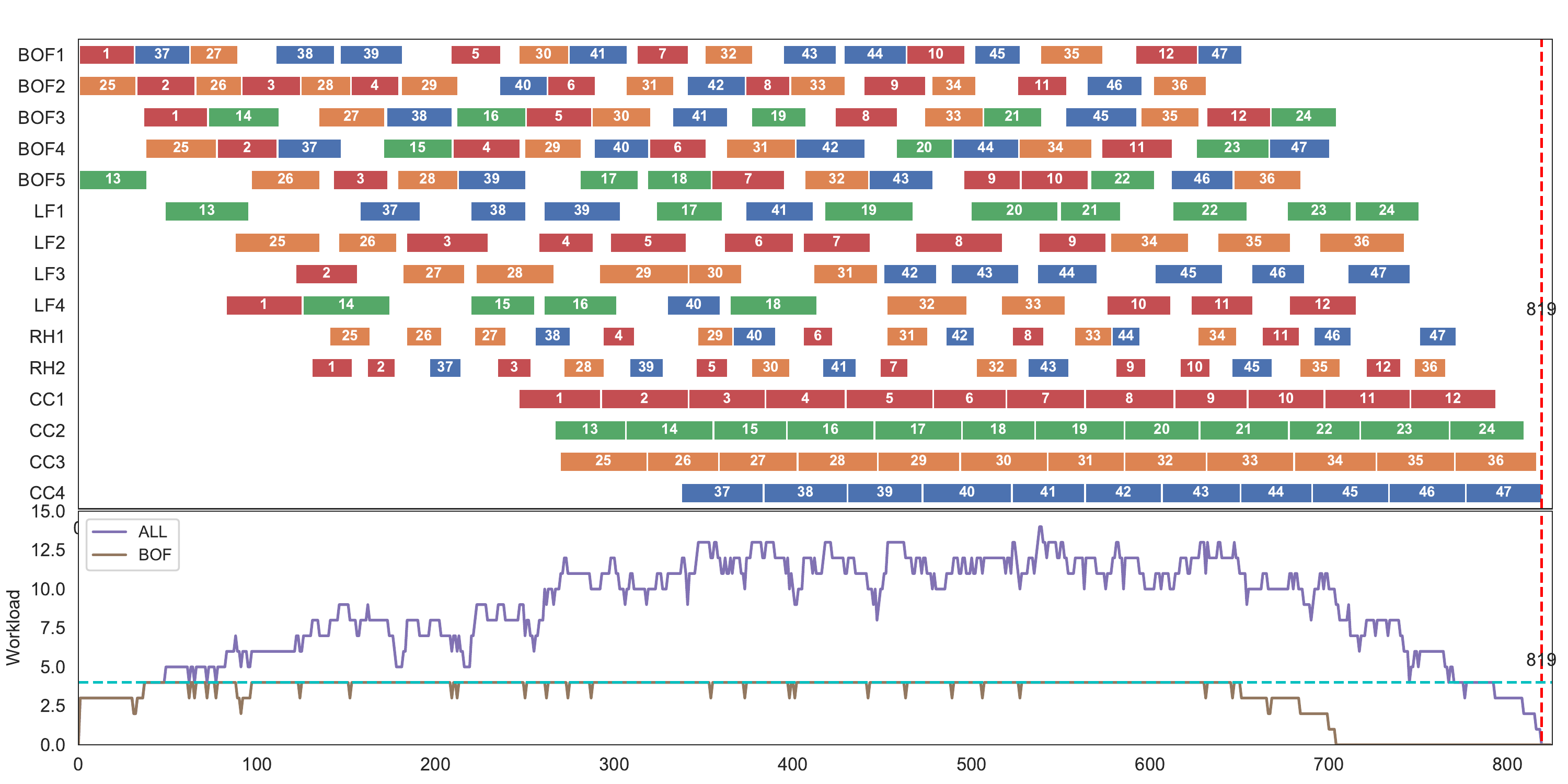}
		\end{minipage}%
	}%
	\caption{Gantt chart of Instance A}
	\label{fig:GanttA}
\end{figure}

\begin{figure}[htbp]
	\centering
	\subfigure[Scenario 1: $C_{\max}=697,W_{\rm{tot}}=60$]{
		\begin{minipage}[t]{1.0\linewidth}	
			\centering		
			\includegraphics[width=1.0\linewidth]{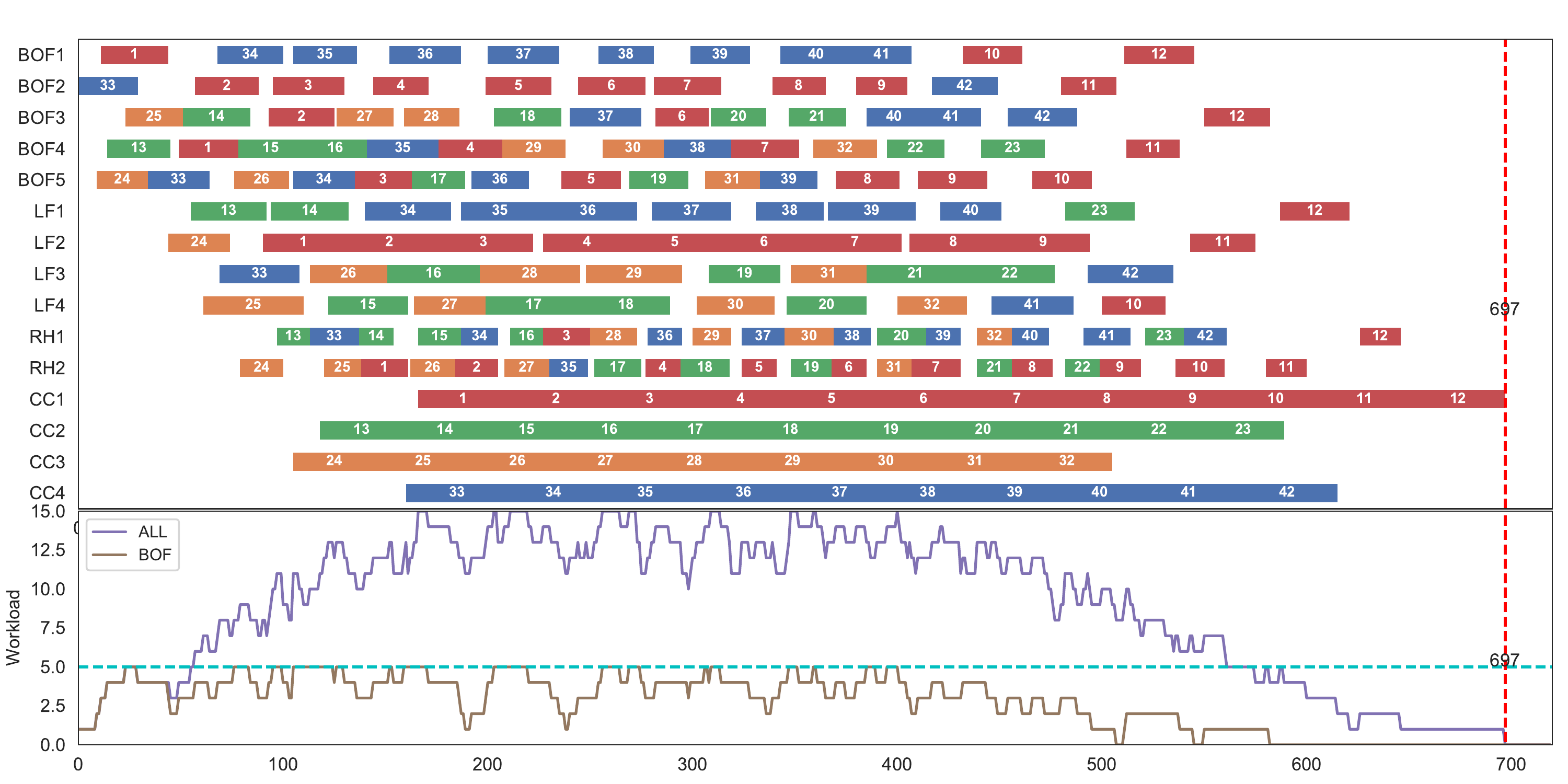}			
		\end{minipage}%
	}\\	
	\subfigure[Scenario 2: $C_{\max}=741,W_{\rm{tot}}=1551$]{
		\begin{minipage}[t]{1.0\linewidth}	
			\centering		
			\includegraphics[width=1.0\linewidth]{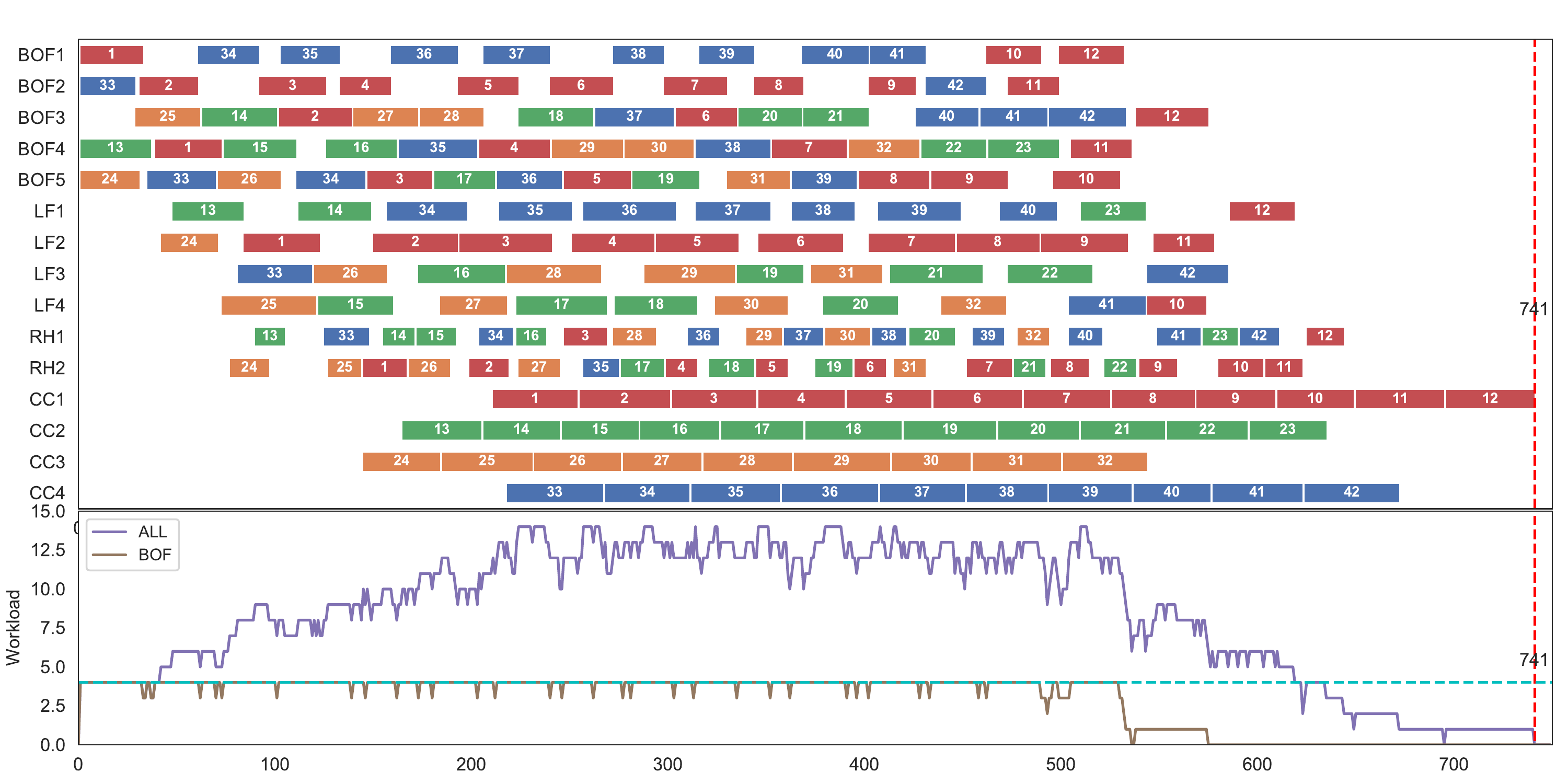}
		\end{minipage}%
	}%
	\caption{Gantt chart of Instance B }
	\label{fig:GanttB}
\end{figure}

\begin{figure}[htbp]
	\centering
	\subfigure[Scenario 1: $C_{\max}=723,W_{\rm{tot}}=455$]{
		\begin{minipage}[t]{1.0\linewidth}	
			\centering		
			\includegraphics[width=1.0\linewidth]{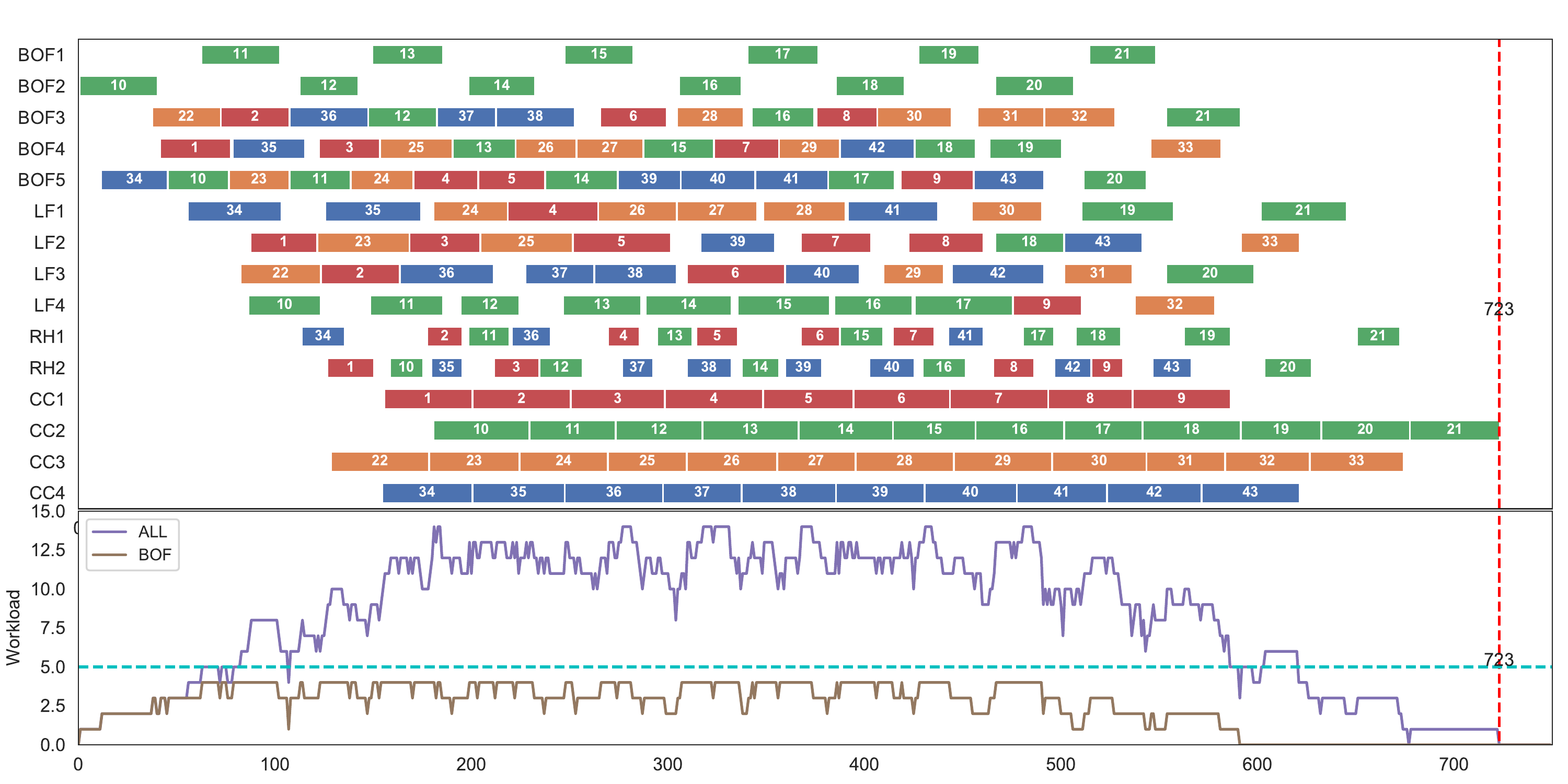}
		\end{minipage}%
	}\\		
	\subfigure[Scenario 2: $C_{\max}=729,W_{\rm{tot}}=2841$]{
		\begin{minipage}[t]{1.0\linewidth}	
			\centering		
			\includegraphics[width=1.0\linewidth]{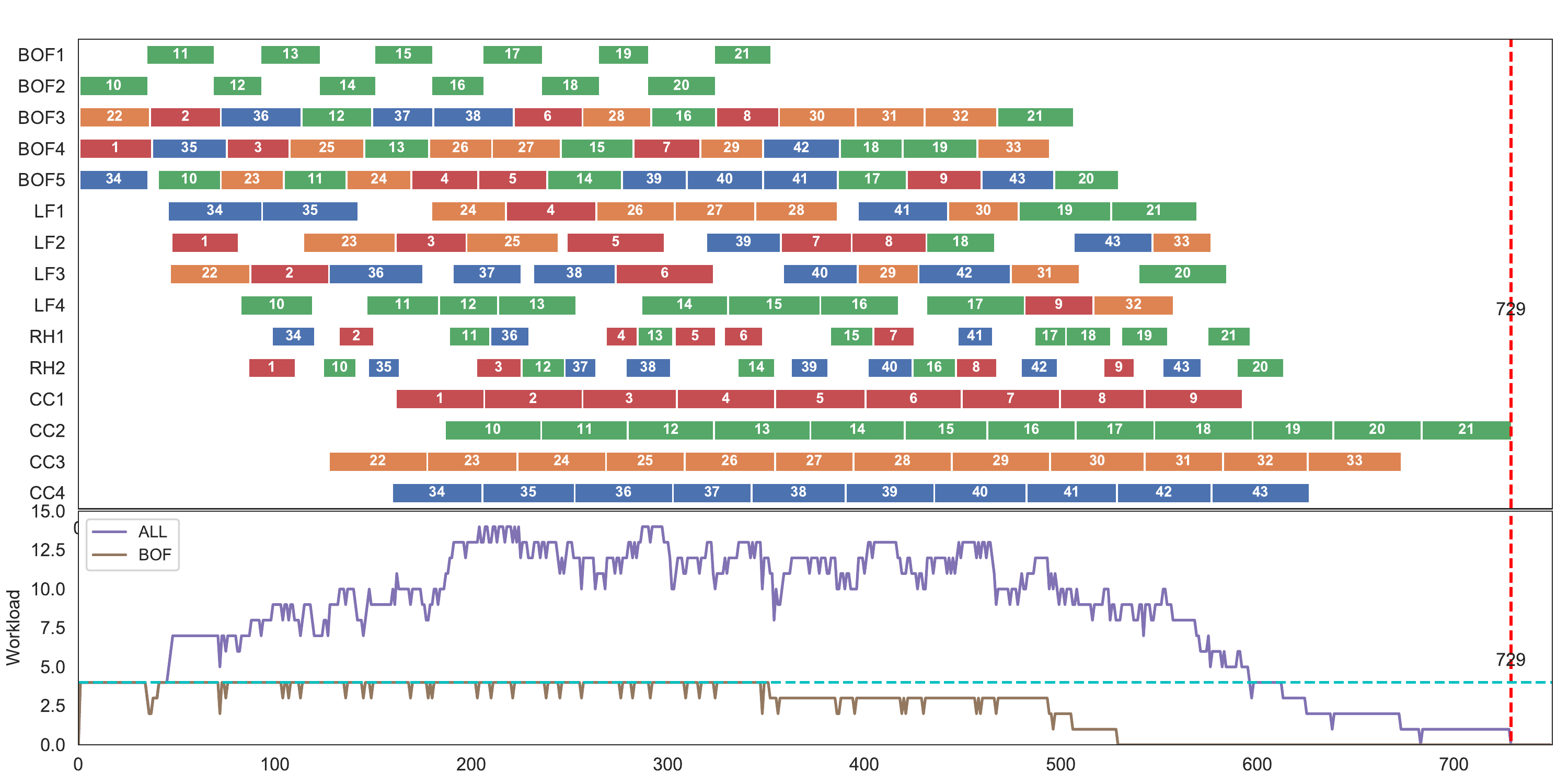}
		\end{minipage}%
	}%
	\caption{Gantt chart of Instance C}
	\label{fig:GanttC}
\end{figure}

\figref{fig:GanttA} -\figref{fig:GanttC} show that the Gantt charts under each scenario of instance A-C. The results indicate optimal scheduling objectives ($C_{\max}$ and $W_{\rm{tot}}$) deteriorate when the capacity is constrained. To observe the impacts on energy demands, we randomly synthesized oxygen consumption of OUS with 32 time periods each of which has 15 minutes and lists their statistical information in Table \ref{tab:inst}. Each instance represents a different level of oxygen demands.

\begin{table}
	\caption{Setting of $d_{\theta}$ for generating test instances}
	\begin{tabular*}{8cm}{llll}
		\hline
		Instance & mean & std. & level \\
		\hline
		A  & 3.98$\times10^4$ & 5.08$\times10^3$ & high  \\
		B  & 3.50$\times10^4$ & 3.46$\times10^3$ & medium  \\
		C  & 3.15$\times10^4$ & 3.32$\times10^3$ & low  \\
		\hline
	\end{tabular*}
	\label{tab:inst}
\end{table}

\begin{figure}[htbp]
	\centering	
	\subfigure[Scenario 1]{
		\begin{minipage}[t]{0.80\linewidth}
			\centering
			\includegraphics[width=1.0\linewidth]{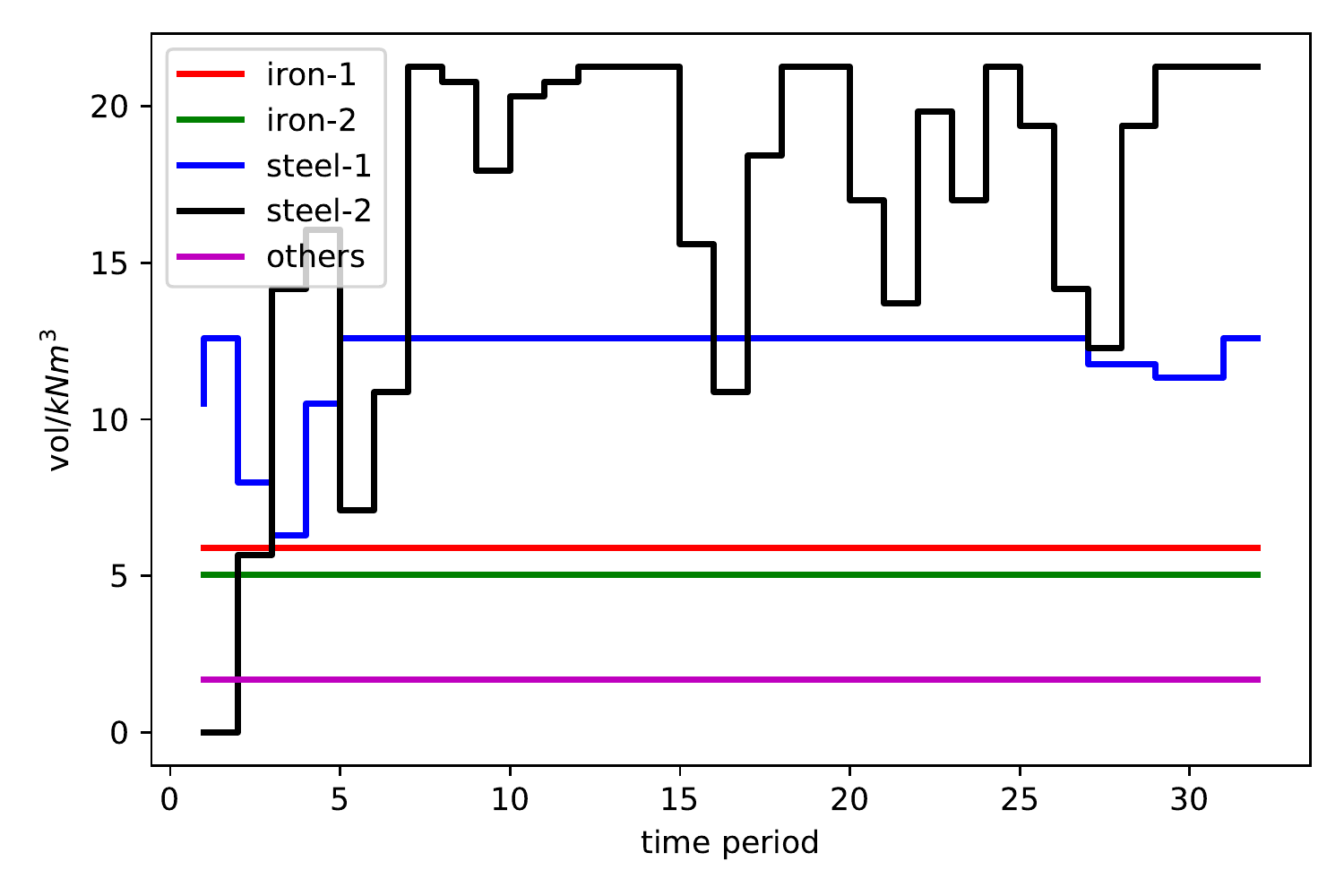}
		\end{minipage}
	}\\
	\subfigure[Scenario 2]{
		\begin{minipage}[t]{0.80\linewidth}
			\centering
			\includegraphics[width=1.0\linewidth]{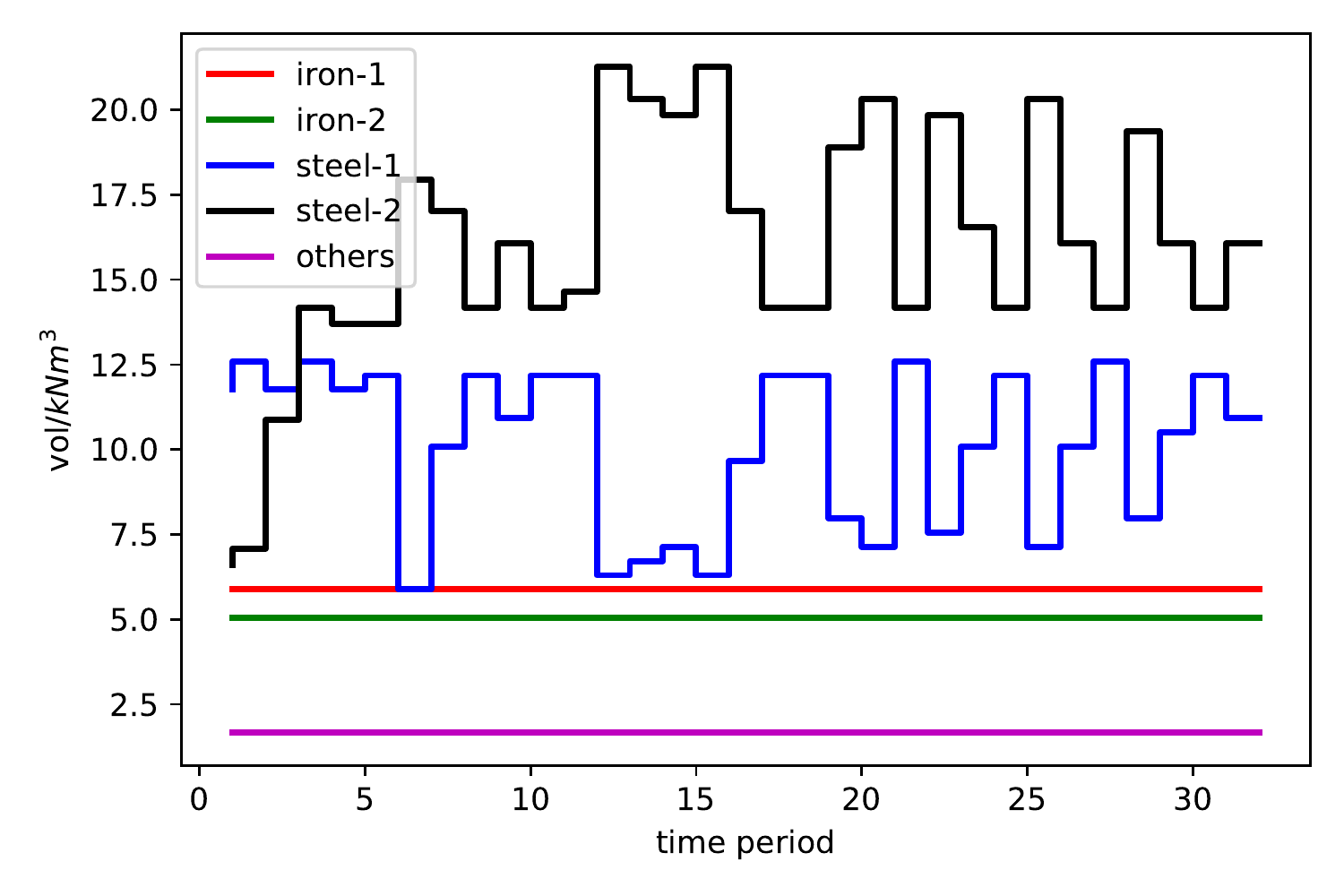}
		\end{minipage}
	}%
	\centering
	\caption{Oxygen demand curve of Instance A}
	\label{fig:DmdA}
\end{figure}

\begin{figure}[htbp]
	\centering	
	\subfigure[Scenario 1]{
		\begin{minipage}[t]{0.80\linewidth}
			\centering
			\includegraphics[width=1.0\linewidth]{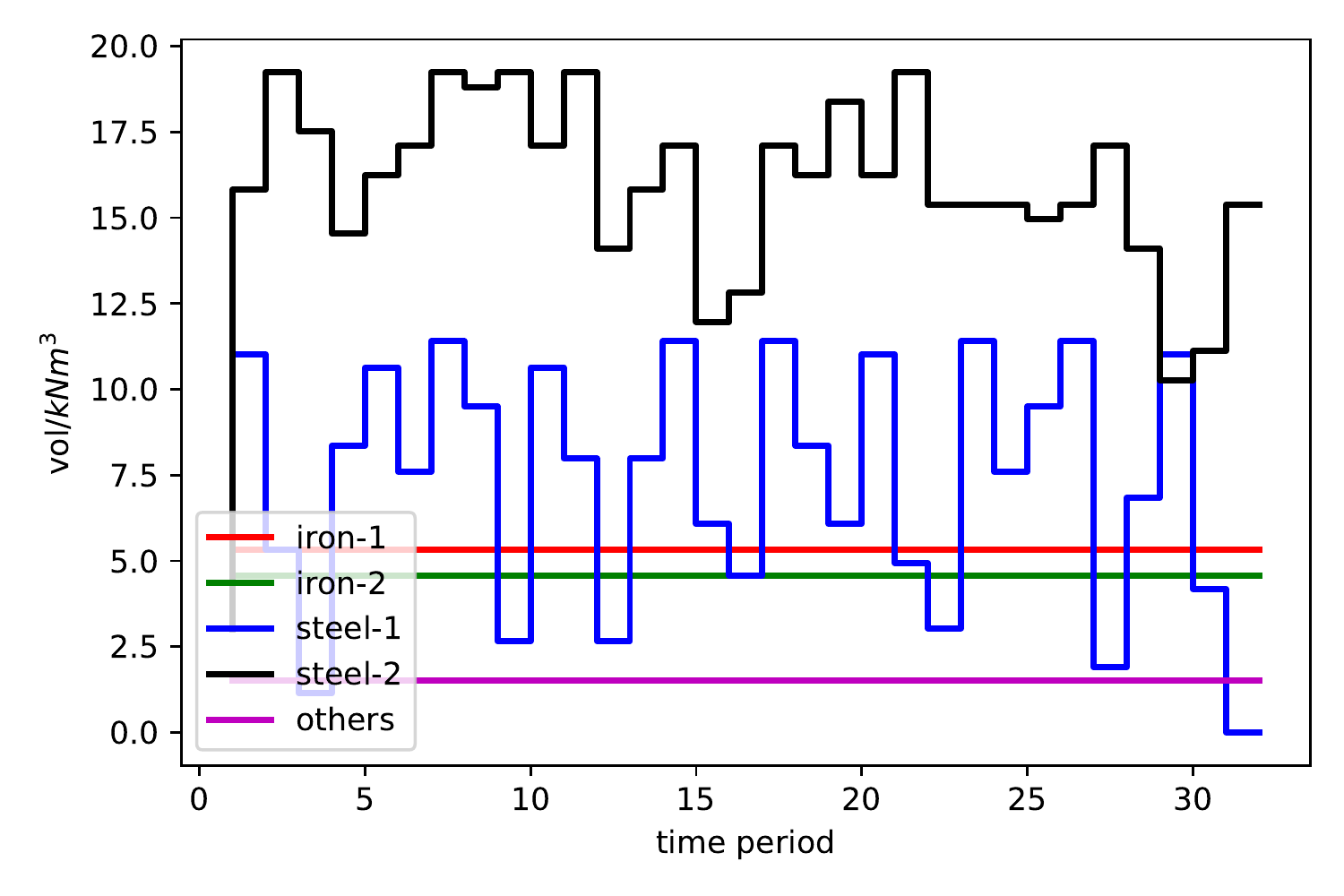}
		\end{minipage}
	}\\
	\subfigure[Scenario 2]{
		\begin{minipage}[t]{0.80\linewidth}
			\centering
			\includegraphics[width=1.0\linewidth]{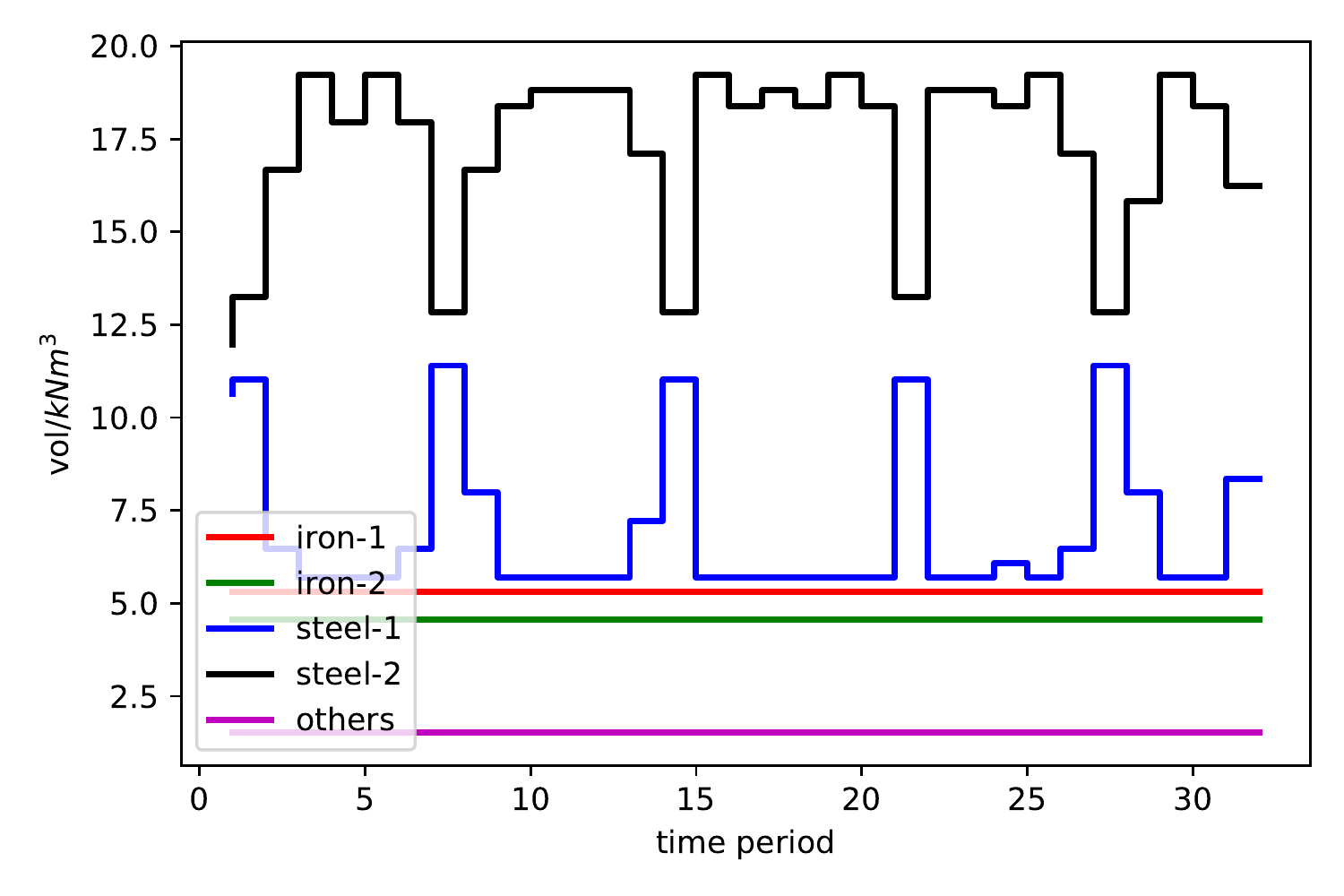}
		\end{minipage}
	}%
	\centering
	\caption{Oxygen demand curve of Instance B}
	\label{fig:DmdB}
\end{figure}

\begin{figure}[htbp]
	\centering	
	\subfigure[Scenario 1]{
		\begin{minipage}[t]{0.80\linewidth}
			\centering
			\includegraphics[width=1.0\linewidth]{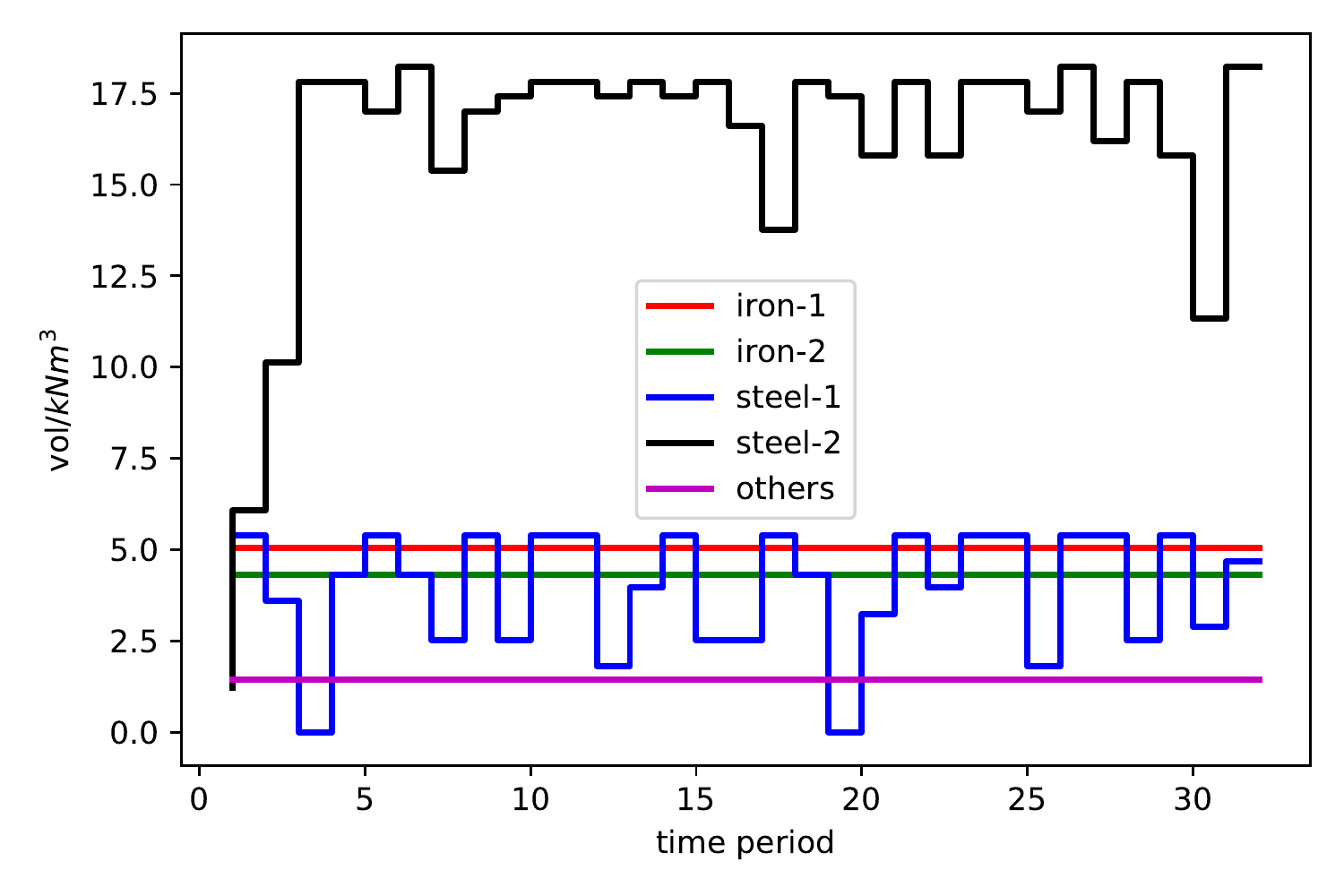}
		\end{minipage}
	}\\
	\subfigure[Scenario 2]{
		\begin{minipage}[t]{0.80\linewidth}
			\centering
			\includegraphics[width=1.0\linewidth]{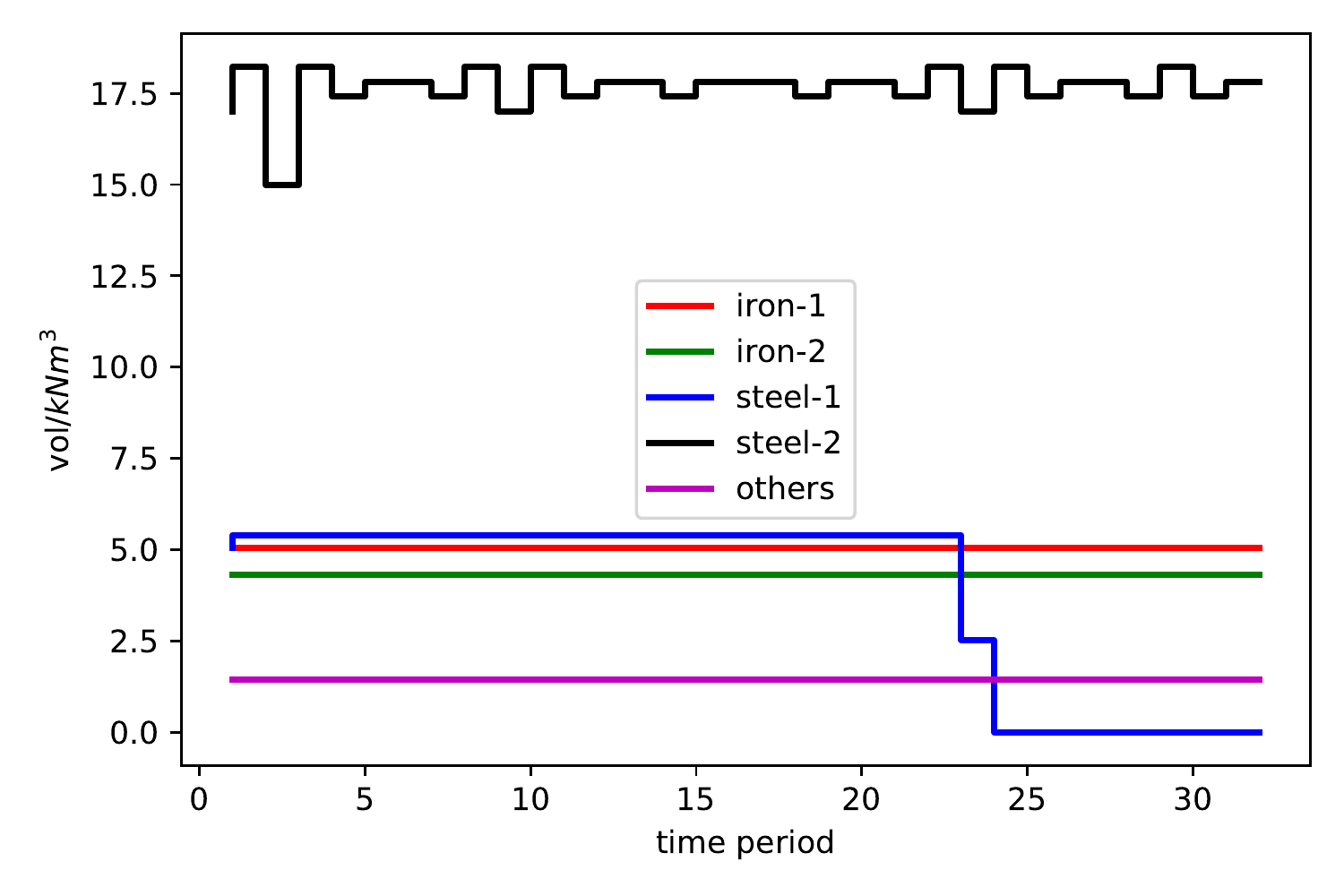}
		\end{minipage}
	}%
	\centering
	\caption{Oxygen demand curve of Instance C}
	\label{fig:DmdC}
\end{figure}

We plotted the demand curves of each user based on deterministic estimates. \figref{fig:DmdA}-\figref{fig:DmdC} show the demand curves under each steelmaking production scheduling scenario of instance A-C, and peak demands not overlapped when the maximum workload is lower. Therefore, the proposed capacity-constrained rescheduling method can effectively shift the demand peak and provide multi-scenario demands for the oxygen distribution model.  

\subsection{Analysis of TSRO model}
To analyse the performance of the TSRO model, we extracted model parameters from the database of the EMS and listed them as follows:
\begin{itemize}
	\item Load of each ASU: $\omega_{r}^{\min} = 1.5 \times 10^4 \rm{Nm^3},\omega_{r}^{\max} =2.0 \times 10^4 \rm{Nm^3} $.
	\item Maximum ramp: $\hat{\omega}^{\max}= 300 \rm{Nm^3}$.
	\item Adjustment rate of continuous users: $\rho_{q}^{\min}=0.8, \rho_{q}^{\max}$ $=1.2$.
	\item Scheduling scenario of discrete users: $S = \left\lbrace 1, 2 \right\rbrace $
	\item Total capacity of OSS: $GV = 6 \times 10^4 \rm{Nm^3}$.
	\item Initial volume of OSS: $GV_{\min}=0.1GV , GV_{\max}= 0.9GV$, $GV_{\rm{mid}}= 0.5GV$
	\item Objective coefficients: $\gamma_1=1.0,\gamma_2=2.0,\gamma_3=20.0$
	\item Number of the time period: $\left| \Theta \right|  = 32$.	
\end{itemize}

In following studies, we assumed the estimated oxygen demands of each user (as shown in \figref{fig:DmdA}-\figref{fig:DmdC}) are randomly varied between $[\bar{d}_{\theta} - \hat{d}_{\theta}, \bar{d}_{\theta} + \hat{d}_{\theta}]$ and $\hat{d}_{\theta}=\eta \times \bar{d}_{\theta} $, where $\eta$ has two candidate levels: $\eta_1 = 0.05, \eta_2 = 0.08$. The two levels respectively cover $68\%, 90\%$ confident intervals.

\subsubsection{Choice of the budget levels}

As discussed in Section \ref{S:RO2D}, we have introduced budget parameters $\Gamma_{\theta}$ and $\Gamma_{\max}$ in TSRO to control the trade-off between the probability of constraint violation and the optimality of the objective function. We define the time-varied budget as follows:
$$\Gamma_{\theta}= \min \left\lbrace Z_{1-a}\sqrt{\theta}+1,  b \left| \Theta \right|  \right\rbrace  $$
where $a$ denotes risk level, $Z$ is its quantile of the standard normal distribution, and $b$ is a constant coefficient. The first term is used to control robustness, and the second term is set for the upper bound of $\Gamma_{\theta}$ to keep feasibility.

\begin{figure}[htbp]
	\centering
	\subfigure[$\eta=0.05$, critical point $(a = 0.05, b=0.40)$]{
		\begin{minipage}[t]{1.0\linewidth}			
			\includegraphics[width=0.96\textwidth]{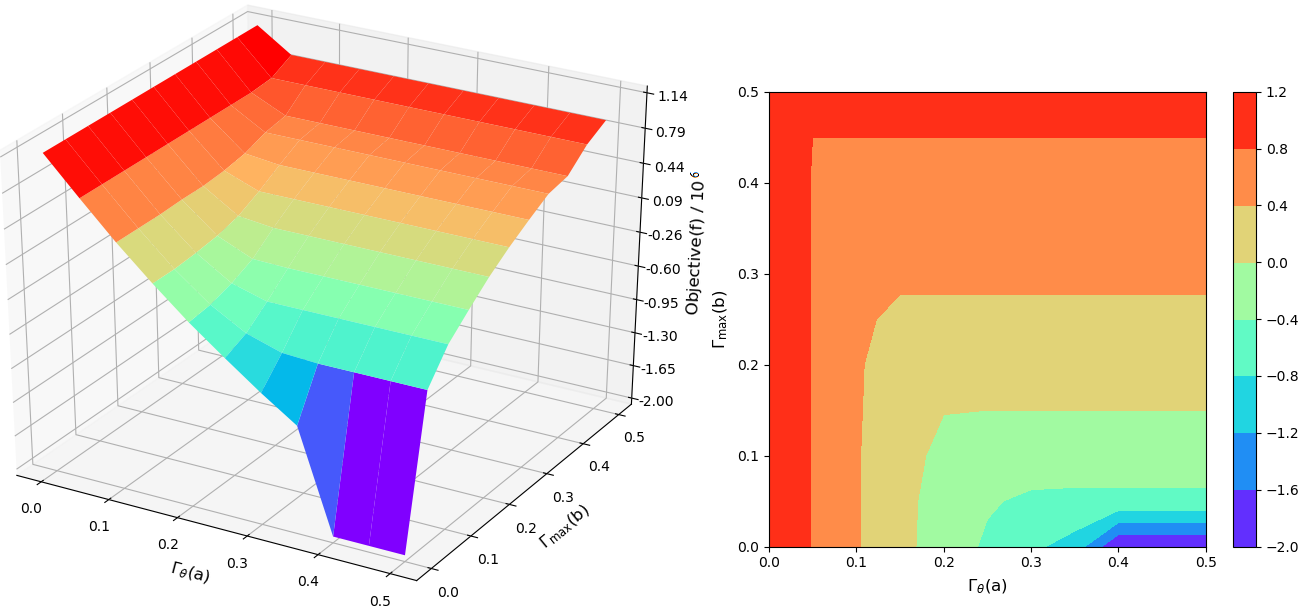}
		\end{minipage}%
	}\\
	\centering
	\subfigure[$\eta=0.08$, critical point $(a = 0.10, b=0.25)$]{
		\begin{minipage}[t]{1.0\linewidth}			
			\includegraphics[width=0.96\textwidth]{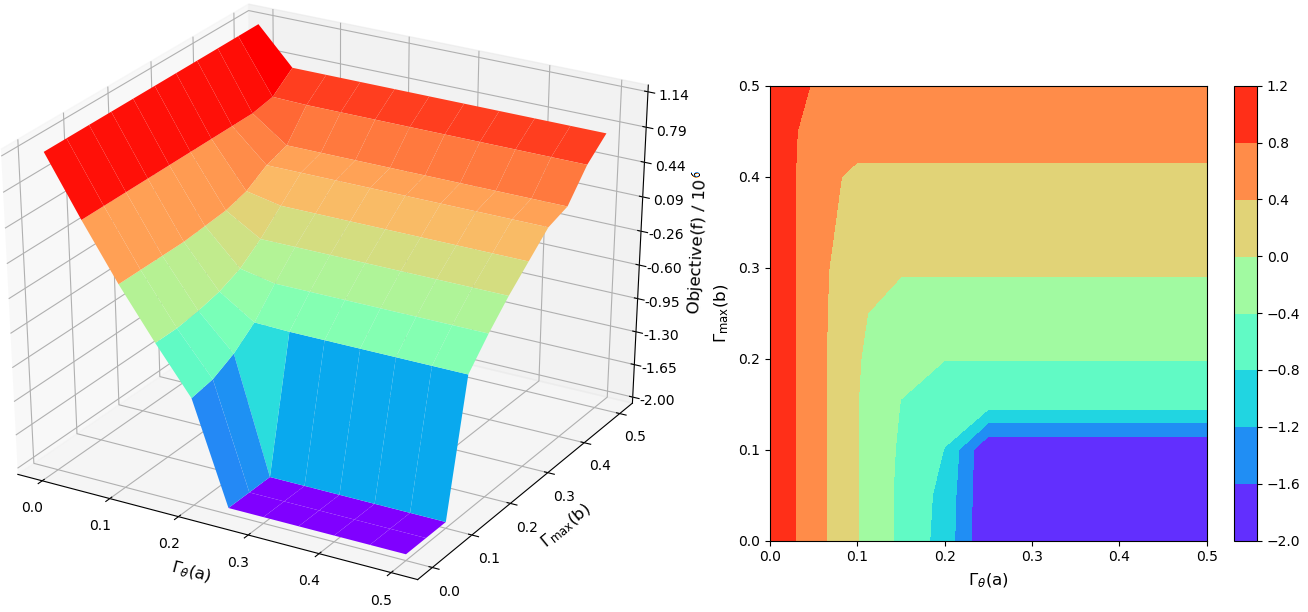}
		\end{minipage}%
	}\\
	\caption{Response under different budget levels}
	\label{fig:risk}
\end{figure}

To choose a appropriate combination of $(a,b)$, we adopted instance B as the test representative and selected $a$ and $b$  from $\left[ 0, 0.5\right]$ with step 0.05, hence get $11 \times 11 =121$ candidate budgets. Then, We ran TSRO with the different value of  $(a,b)$ and plotted the response surfaces of two deviations ($\eta_1=0.05$ and $\eta_2=0.08$), as shown in \figref{fig:risk}. If TSRO is infeasible, let $f=-2.0\times 10^6$. These results demonstrate the objective value of $f$ increases as the risk level decrease, as the maximum budget decreases. When $a \geq 0.05$ or $b \leq 0.40 $, no solution under deviation of $0.05$ is infeasible. When $a \geq 0.10$ or $b \leq 0.25 $, no solution under deviation of $0.08$ is infeasible. According these critical points shown in the figures, we recommend $a=0.10, b= 0.40$ in following computational studies.
 
\subsubsection{Sensitivity analysis on ramp capacity}
In this section, we study the relationship between the ramp capacity of ASU and the performance of robust models. The intuition is that higher ramp capacity let the EMS have a better ability to deal with a high variation of oxygen. Therefore, we attempted to verify how much benefit the TSRO model provides under different ramp capacities between $\left[0, 500 \right] Nm^3$ with the step of $50 Nm^3$. 

\begin{figure}[htbp]
	\centering
	\subfigure[$\eta_1 = 0.05$]{
		\begin{minipage}[t]{1.0\linewidth}			
			\includegraphics[width=0.96\textwidth]{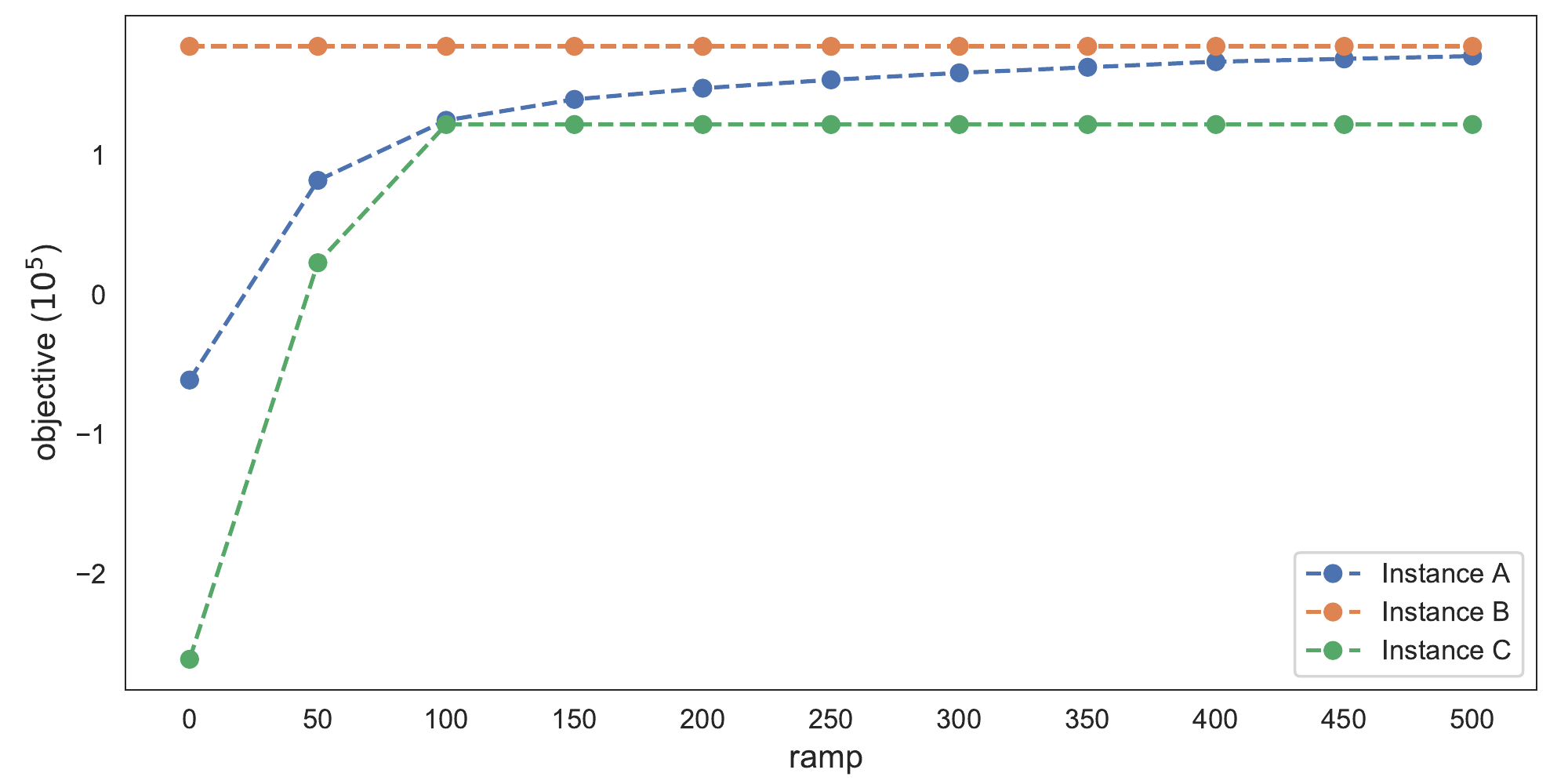}
		\end{minipage}%
	}\\
	\centering
	\subfigure[$\eta_2 = 0.08$]{
		\begin{minipage}[t]{1.0\linewidth}			
			\includegraphics[width=0.96\textwidth]{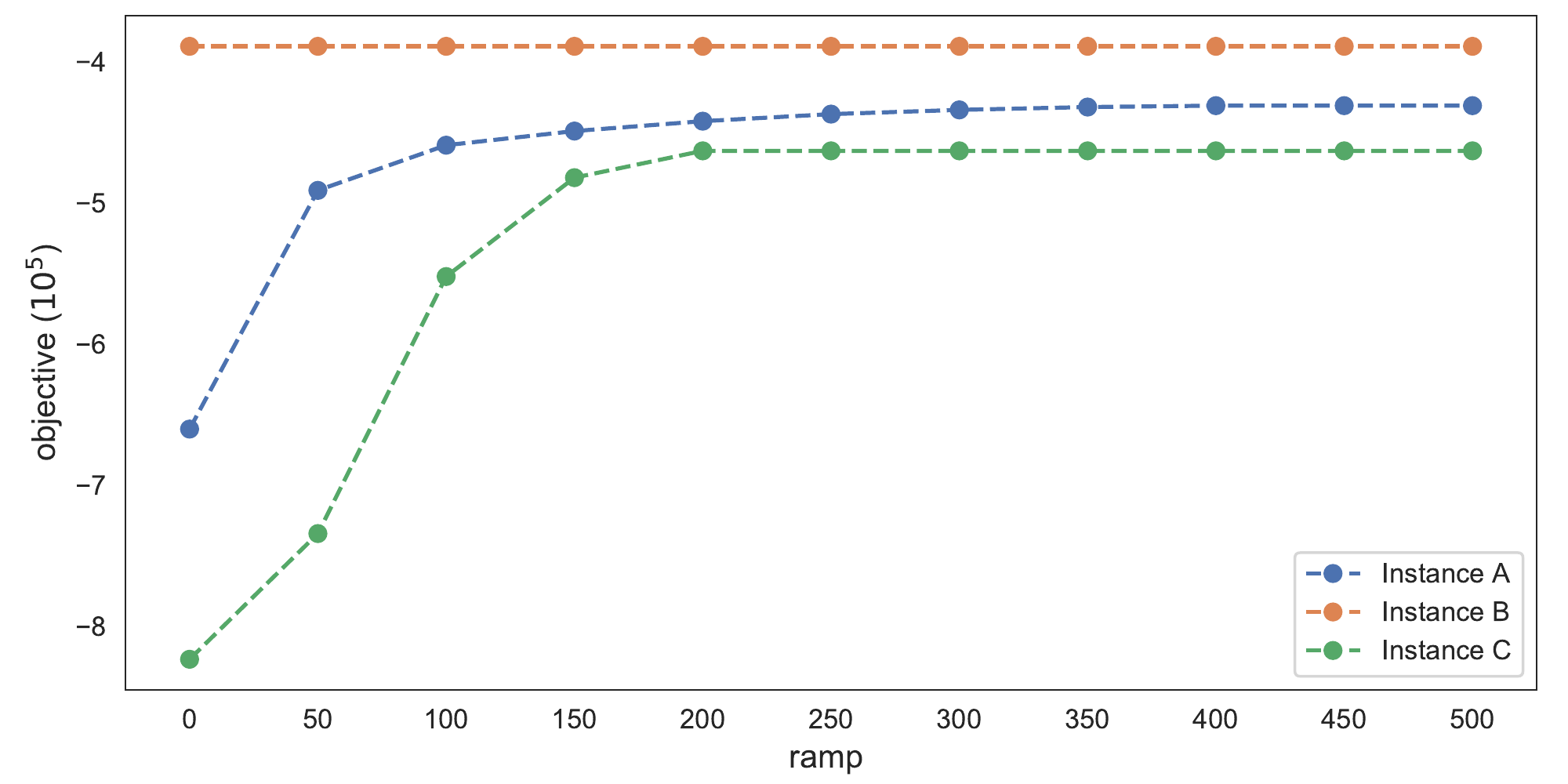}
		\end{minipage}%
	}\\
	\caption{Response under different ramp capacities}
	\label{fig:ramp}
\end{figure}

\figref{fig:ramp} summarizes the computational results for two groups studies: the deviation level of $\eta_1=0.05$, and $\eta_2=0.08$. The results indicate that the TSRO models for instance A and C can obtain more optimal objectives by enlarging the ramping capacity when ramping capacity $\hat{\omega}_{\max} \leq 200$, but the model for instance B keeps its objective unchanged. It demonstrates clearly that the model with higher and lower demands needs to design an appropriate ramp capacity to obtain more profits.  

\subsubsection{Robustness comparison by simulation}

With the budget-based uncertainty set $\mathcal{U}$, TSRO model provides a protective solution. We randomly sampled two demand vectors $d^{1}_{\theta}$ and $d^{2}_{\theta}$ following a truncated normal distribution, which fell between $\left[\bar{d}_{\theta} -\hat{d}_{\theta}, \bar{d}_{\theta} + \hat{d}_{\theta}\right] $, and let the random demand sample $d_{\theta}=(d^{1}_{\theta}+d^{2}_{\theta})/2$ . Notice that, the random sample may locate outside of the budget-based uncertainty set. In this way, we can test the robustness of the proposed TSRO model. In following studies, we test three instances (A-C) with five different initial gasholder levels ($GV_0 = (0.3 \sim 0.7) GV $),  hence get $3\times5 = 15$ cases under each deviation ($\eta \in \left\lbrace 0.05, 0.08 \right\rbrace$). Then, we independently run the distribution model with deterministic optimization (DO), the proposed TSRO with the fixed budget level ($a=0.10, b=0.40$), and stochastically simulate the solutions of DO and TSRO  with 1000 rounds (-simulated). 

\begin{figure}[htbp]
	\centering
	\subfigure[$\eta_1 = 0.05$]{
		\begin{minipage}[t]{1.0\linewidth}			
			\includegraphics[width=0.96\textwidth]{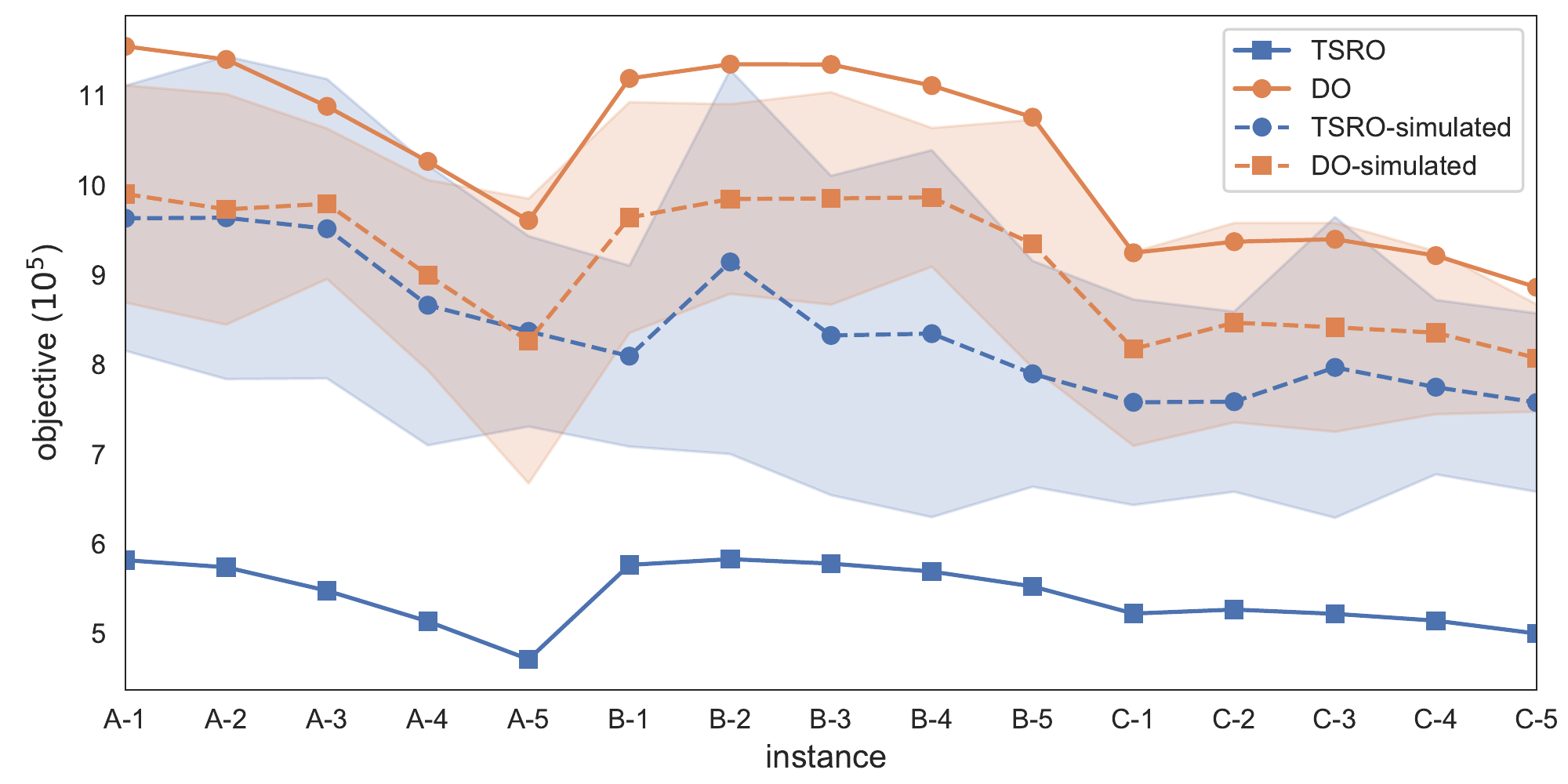}
		\end{minipage}%
	}\\
	\centering
	\subfigure[$\eta_2 = 0.08$]{
		\begin{minipage}[t]{1.0\linewidth}			
			\includegraphics[width=0.96\textwidth]{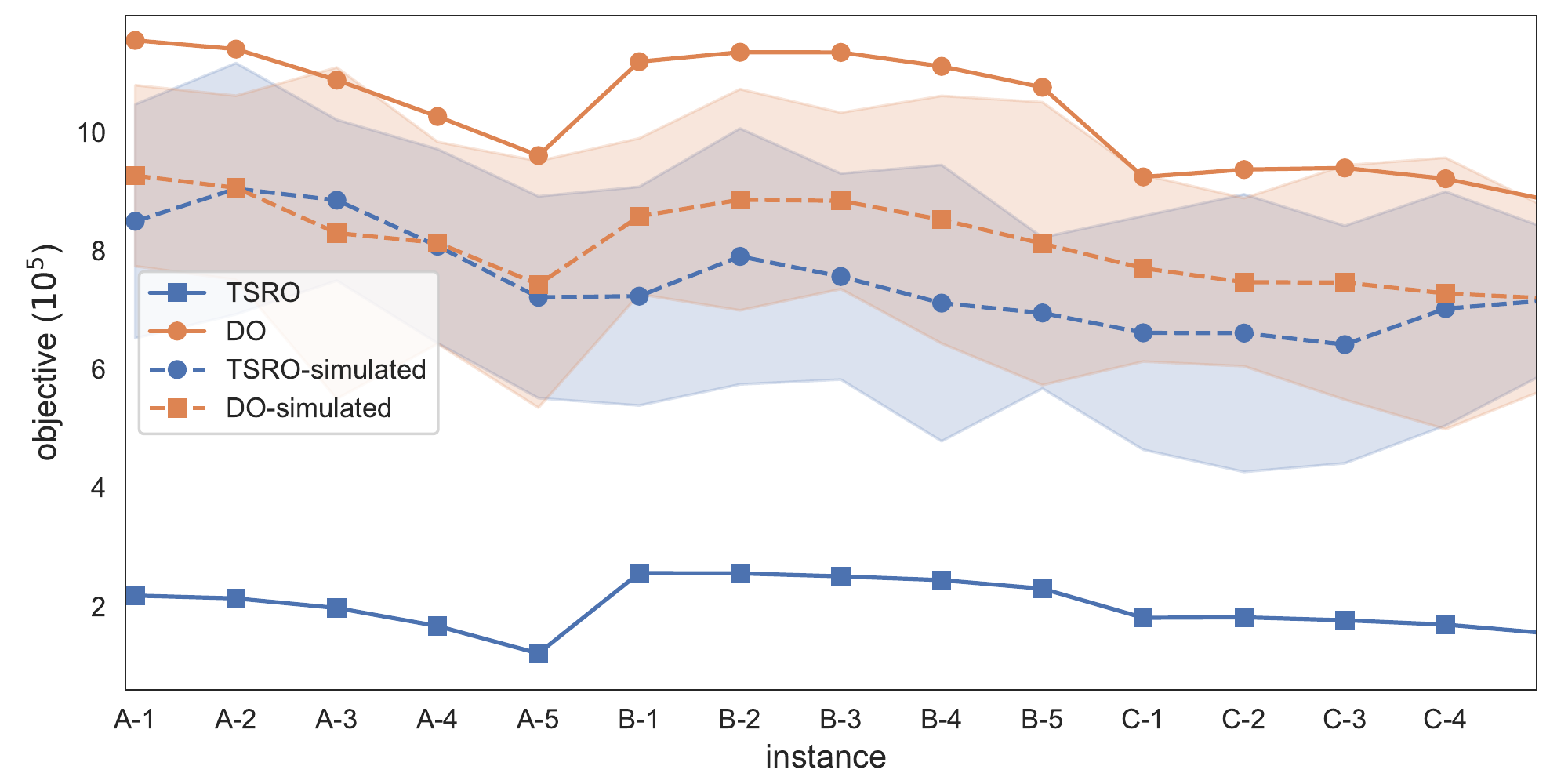}
		\end{minipage}%
	}\\
	\caption{Simulation results of different cases}
	\label{fig:robust}
\end{figure}

\figref{fig:robust} plots the objective values of DO and RO, and their average objectives under simulation with 95\% confidence intervals (2 standard deviations). The results indicate the objective of TSRO is lower than DO, but it provides promising protection that is lower than the 95\% confidence intervals. Under stochastic environments, the simulated objective of TSRO is only slightly lower than simulated DO, which indicate it less conservative. Not surprisingly, the deterministic optimization fails in hedging against uncertainty for most instances except in some cases of instance B (with medium demand). This means that too big or too small demand will increase the risk of oxygen distribution.

\section{Conclusions}
\label{S:7}
In this paper, we studied an optimal oxygen distribution problem with flexible demands under uncertainty and proposed a TSRO model with a budget-based uncertainty set that allows protecting the initial distribution policy with low conservatism. With computational studies, the following conclusions can be drawn:
\begin{enumerate}[(1)]
\item The GP-based time series model can forecast varied oxygen demands with 95\% confidence intervals. The capacity-constrained rescheduling model can effectively shift the step-wise oxygen demands. 	
\item According to the level of risk aversion, the decision-maker can select appropriate budget parameters for the TSRO model. The ramp capacity also contributes to the flexibility of oxygen distribution.
\item Comparing deterministic optimization commonly used in existing studies, the TSRO model can provide a more safe solution, but its average performance under uncertainty has not deteriorated.
\end{enumerate}
  
In conclusion, model simulation is coherent with the behaviour of our design and confirm that the proposed model provides a promising way to solve the uncertain oxygen distribution problem. In the future, the following issues will remain to be further studied:
\begin{enumerate}[(1)]
\item In the investigated problem, we only consider as uncertain energy demands from the manufacturing system. However, many other uncertain factors should be further studied, such as energy market demands and machine breakdown.
\item In the proposed TSRO model, we only know the interval of the uncertain demands, but the probabilistic information has not embedded the model. So, how to build a data-driven uncertainty set is another challenging issue.
\end{enumerate}

\section*{Acknowledgment}

This research was financially supported by the by the National Natural Science Foundation of China under Grant 61873042, in part by the Open Fund of Guangxi Key Laboratory of Automatic Detection Technology and Instrument (No. YQ19202). Thanks to the China Scholarship Council (CSC) for supporting a 1-year visit at University College London (UCL), U.K.

\nomenclature[D01]{$R,r$}{Suppliers of OGS.}
\nomenclature[D02]{$Q,q$}{Users of OUS,and $Q^N, Q^E$ repectively denote continuous and discrete user sets.}
\nomenclature[D03]{$\Theta,\theta$}{Time periods of oxygen distribution, $\Theta_{\theta}={1,\dots,\theta}$.}
\nomenclature[D04]{$S, s$}{Scenario set of discrete engery demands.}
\nomenclature[D05]{$G, g$}{Stage set and index of a steelmaking work.}
\nomenclature[D06]{$M_g,m$}{Machine set and index in stage $g$.}
\nomenclature[D07]{$J,j$}{Job set relased into a steelmaking work.}
\nomenclature[D08]{$A$}{Task set processed in the steelmaking workshop.}
\nomenclature[D09]{$H,h$}{Bathes groupped by the released jobs.}
\nomenclature[D10]{$L,l $}{Time units of the process scheduling horizon.}
\nomenclature[D11]{$T,t $}{Time points of the energy demand time series.}

\nomenclature[P01]{$\omega_{r}^{\min}$}{Minimum load of oxygen generator $r$, $\left[ Nm^3\right] $}
\nomenclature[P02]{$\omega_{r}^{\max}$}{Maximum load of oxygen generator $r$, $\left[ Nm^3\right] $}
\nomenclature[P03]{$\hat{\omega}_{\max}$}{Maximum deviation of oxygen load between period $(\theta-1)$ and $\theta$, $\left[ Nm^3\right] $}		
\nomenclature[P04]{$\rho_{q}^{\min}$}{Minimum adjustment rate of oxygen user $q$, $\left[ \%\right] $}	
\nomenclature[P05]{$\rho_{q}^{\max}$}{Maximum adjustment rate of oxygen user $q$, $\left[ \%\right] $}	
\nomenclature[P06]{$GV$}{Total capacity of the oxygen storage, $\left[ Nm^3\right] $}
\nomenclature[P07]{$GV_{\min}$}{Minimum level of the oxygen storage, $\left[ Nm^3\right] $}
\nomenclature[P08]{$GV_{\max}$}{Maximum level  of the oxygen storage, $\left[ Nm^3\right] $}
\nomenclature[P09]{$GV_{\rm{mid}}$}{Middle level of the oxygen storage, $\left[ Nm^3\right] $}
\nomenclature[P10]{$d_{q,\theta}$}{Demand of energy user $q$ at period $t$,  $\left[ Nm^3\right] $}
\nomenclature[P11]{$\Gamma_{\theta},\Gamma_{\max}$}{Time-variant budget and maximum budget of uncertainty set}
\nomenclature[P12]{$PT_{g,j}$}{Processing time of task $A_{g,j}$, $\left[ min \right]$}
\nomenclature[P13]{$TT_{g,gg}$}{Transfer time between stage $g$ and $gg$,$\left[ min \right] $}
\nomenclature[P14]{$SU_{h}$}{Setup time of batch $h$, $\left[ min \right] $}
\nomenclature[P15]{$ES_{g,j}$}{Earliest starting time point of task $A_{g,j}$}
\nomenclature[P16]{$LF_{g,j}$}{Latest finish time point of task $A_{g,j}$}
\nomenclature[P17]{$WL_{q}^{\max}$}{Maximum capacity of oxygen in the $q^{th}$ scenario of the primary steelmaking stage}
\nomenclature[P18]{$TL$}{Time lags of state space model in time series}
\nomenclature[P19]{$N$}{Size of training input in time series}
\nomenclature[P20]{$k, K$}{Kernel function of the GP-based forecasting model.}

\nomenclature[V1]{$\omega_{r,\theta}$}{Load of gaseous oxygen of supplier $r$ at time period $\theta$, $Nm^3$}
\nomenclature[V2]{$\rho_{q,\theta}$}{Oxygen flow rate of of user $q$ at time period $\theta$}
\nomenclature[V3]{$\delta_{\theta}$}{Deviation from  the gasholder level to its middle level at time period $\theta$, $Nm^3$}
\nomenclature[V4]{$\epsilon_{\theta}$}{Volume of surplus or shortage gas at time period $\theta$, $Nm^3$}
\nomenclature[V5]{$GV_{\theta}$}{Volume of gas storage at time period $\theta$, $Nm^3$}
\nomenclature[V6]{$z_{s}$}{Binary variable which is equal to 1 if and only if the steelmaking schedule of the $s^{th}$ scenario is selected}
\nomenclature[V7]{$x_{g,j,l}$}{Binary variable which is equal to 1 if and only if task $A_{g,j}$ is in-process at time $t$}
\nomenclature[V8]{$C_{g,j}$}{Completion time of task $A_{g,j}$}

\printnomenclature

\bibliography{mybibfile}
\end{document}